\documentclass[letterpaper,12pt]{amsart}
\usepackage[dvipsnames]{xcolor}
\usepackage{latexsym,array,delarray,amsthm,amssymb,epsfig,setspace,tikz,amsmath,enumerate,mathrsfs,graphicx, mathtools}
\usepackage[shortlabels]{enumitem}
\usepackage{ytableau}
\usepackage[ruled,vlined]{algorithm2e}
\usepackage{blkarray}
\usepackage{etoolbox}

\usepackage{url}
\usepackage{microtype}
\usepackage[colorlinks=true,hyperindex, linkcolor=magenta, pagebackref=false, citecolor=cyan]{hyperref}
\usepackage[alphabetic,lite]{amsrefs}

\usepackage{subcaption}

\usepackage[shortlabels]{enumitem}
\usepackage{arydshln}
\usepackage{cleveref}
\usepackage[backgroundcolor=white,bordercolor=black]{todonotes}
\usepackage{comment}
\usepackage{pst-node}
\usepackage{tikz-cd} 
\usetikzlibrary{arrows,automata}

\theoremstyle{plain}

\newtheorem{thm}{Theorem}[section]
\newtheorem{theorem}[thm]{Theorem}
\newtheorem{lemma}[thm]{Lemma}
\newtheorem{prop}[thm]{Proposition}
\newtheorem{corollary}[thm]{Corollary}
\newtheorem{conjecture}[thm]{Conjecture}
\newtheorem*{thm*}{Theorem}
\newtheorem*{lemma*}{Lemma}
\newtheorem*{prop*}{Proposition}
\newtheorem*{cor*}{Corollary}
\newtheorem*{conj*}{Conjecture}

\theoremstyle{definition}

\newtheorem{definition}[thm]{Definition}
\newtheorem{example}[thm]{Example}

\newtheorem{hyp}[thm]{Hypothesis}

\theoremstyle{remark}

\newtheoremstyle{case}{}{}{}{}{}{:}{ }{}
\theoremstyle{case}

\usepackage{ upgreek }

\newcommand{\bfx}{\mathbf{x}}

\newcommand{\cG}{\mathcal{G}}
\newcommand{\Gcal}{\mathcal{G}}

\newcommand{\Lcal}{\mathcal{L}}
\renewcommand{\S}{\mathbb{S}}
\newcommand{\Mcal}{\mathcal{M}}

\newcommand{\ind}{\mbox{$\perp \kern-5.5pt \perp$}}

\newcommand{\Gbar}{\overline{\mathcal{G}}}

\newcommand{\lra}{\leftrightarrow}

\newcommand{\Z}{\mathbb{Z}}
\newcommand{\cC}{\mathcal{C}}
\newcommand{\Bcal}{\mathcal{B}}
\newcommand{\Ne}{\mathrm{Ne}_{\Gcal}}

\title[Symmetrically colored graphical models with toric vanishing ideals]{Symmetrically colored Gaussian graphical models with toric vanishing ideals}

\author{Jane Ivy Coons, Aida Maraj, Pratik Misra, Miruna-\c Stefana Sorea}
\begin{document}

\maketitle

\begin{abstract}
\normalsize 
 A colored Gaussian graphical model is a linear concentration model in which equalities among the concentrations are specified by a coloring of an underlying graph. The model is called RCOP if this coloring is given by the edge and vertex orbits of a subgroup of the automorphism group of the graph. 
 We show that RCOP Gaussian graphical models on block graphs are toric in the space of covariance matrices and we describe Markov bases for them.
 To this end, we learn more about the combinatorial structure of these models and their connection to Jordan algebras. 
\end{abstract}

\footnote{{\small \textbf{\emph{Keywords}}: colored Gaussian graphical models, linear concentration models, block graphs, RCOP models, Jordan algebras, toric varieties, Markov bases.}}

\section{Introduction}

Gaussian graphical models are  statistical models in which the relationships between random variables are encoded by a graph. Given a graph $G$ with vertex set $[n] := \{1,\dots,n\}$, the associated Gaussian graphical model is the linear concentration model in which the concentrations corresponding to non-edges of $G$ are constrained to be $0$ and all other concentrations are free. Gaussian graphical models are widely applicable in both the physical and social sciences. For example, they appear in the field of computational biology to model gene interactions \cite{toh2002} and in the field of environmental psychology to model community attitudes towards sustainable behaviors \cite{bhushan2019}. However for a Gaussian graphical model on a graph with many edges, the large number of parameters can make computations intractable. Thus it can be useful to impose further symmetry constraints among the concentrations for covariates that we expect to behave similarly  \cite{uhler2018, wit2015}.

In \cite{hojsgaard2008}, Hojsgaard and Lauritzen introduced several types of graphical models with added symmetries among the parameters. These symmetries are encoded by colorings of the edges and vertices of $G$; we set two entries of the concentration matrix equal to one another if their corresponding edges or vertices have the same color. These models are used to study gene regulatory networks wherein one imposes symmetries among genes with similar expression patterns \cite{toh2002,vinciotti2016}. They have also been applied to the analysis of longitudinal data on the performance of several companies in the same market \cite{abbruzzo2016}. Adding these constraints to a graphical model can reduce the maximum likelihood threshhold of the model; this makes it possible to compute the maximum likelihood estimate with relatively few data points \cite{MRS21,uhler2011}. 

In the present paper, we study the geometry of RCOP models (\Cref{def:RCOPModel}), wherein the coloring of the underlying graph is given by the orbits of a group of automorphisms of $G$ \cite{hojsgaard2008}. We denote the underlying graph along with its coloring by $\Gcal$. In an RCOP model, the concentration matrices $K$ are restricted to lie in a certain linear space $\Lcal_{\Gcal}$. The name ``RCOP" roughly stands for ``restrictions on concentrations generated by permutation symmetries".

We are especially interested in the algebraic variety $\Lcal_{\Gcal}^{-1}$, which is the Zariski closure of the set of all inverses of matrices in $\Lcal_{\Gcal}$. The intersection of $\Lcal_{\Gcal}^{-1}$ with the positive definite cone is the set of all covariance matrices for distributions in the RCOP model. We aim to compute the vanishing ideal of $\Lcal_{\Gcal}^{-1}$, denoted $I_{\Gcal}$, and to determine when it is toric.  Knowledge of the polynomials in this ideal is applicable to several different parameter inference problems. For example, these polynomials are often used to compute the number of complex critical points of the log-likelihood function, known as the ML-degree, of a given model \cite{BCEMR20, dinu2021geometry, manivel2022, SU10}. In \cite{BCEMR20}, the toric structure of the ideal is instrumental to the ML-degree computation.
The generators of this ideal are useful for proving that a candidate function of the data is the maximum likelihood estimate, or MLE. One can show that the candidate MLE lies in the model by checking that it vanishes on the generators of the ideal.
Moreover, the quadratic binomials in this vanishing ideal are utilized in the TETRAD procedure, which is used to infer the structure of the underlying graph given some data \cite{spirtes2000}. Our main result, \Cref{thm:MainSec6}, gives all such quadratic binomials in the case where the RCOP graph is a one-clique sum of complete graphs, or block graph (\Cref{def:block_graph}).

We are now ready to outline the main theorem. Let $\Gcal$ be a connected RCOP block graph with underlying uncolored graph $G$. The vanishing ideal $I_G$ of the set of covariance matrices arising from $G$ is known to be toric and has generating set of quadratic binomials \cite{misra2019gaussian}. We show that the ideal of the RCOP graph $I_\cG$ is generated in a natural way by these quadratics and certain binomial linear forms that are extracted by examining the sequences of colors along paths in $\Gcal$.
Moreover, these linear forms come from some coloring of the complete graph on $n$ vertices, which we denote by $\bar{\Gcal}$ and examine in detail in \Cref{sec:completionsAreRCOP}.
\begin{theorem}\label{thm:IntroMain}
Let $\Gcal$ be an RCOP block graph. Then $I_{\Gcal}$ is toric and generated in degrees one and two. Moreover,
the union of generating sets of $I_G$ and $I_{\bar{\Gcal}}$ 
is
a generating set for~$I_{\Gcal}$.
\end{theorem}
Our proof strategy is as follows: we define a monomial map, called the \emph{shortest path map} on $\cG$, and show that the Zariski closure of its image is equal to that of the rational map defining the model. 
In order to accomplish this, we develop a combinatorial understanding of RCOP block graphs and employ the theory of Markov bases.

\textbf{Structure of the paper.} We start with some preliminaries on   RCOP models and Markov bases in~\Cref{sec:prelim}. In ~\Cref{sec:BlockGraphs} we define block graphs and recall the known Markov basis for the Gaussian graphical model of an uncolored block graph. We end this section with the definition of the shortest path map of an RCOP block graph.  \Cref{sec:completionsAreRCOP} and \Cref{sec:AnotherMonomialMap} contain several key results related to the structure of RCOP block graphs in service of the main result. We introduce the completion of an RCOP block graph in \Cref{sec:Jordan} and show that this completion is itself RCOP. We further show that its vanishing ideal is linear and provide an explicit generating set for it. \Cref{sec:Main} is dedicated to the proof of the main theorem. We conclude the paper with a discussion on future directions.

\section{Preliminaries}\label{sec:prelim}

\subsection{The vanishing ideal of an RCOP  model}
Let $\cG$ be a colored graph on $n$ vertices with edge set $E(\Gcal)$. We assume throughout that $\Gcal$ has no loops or multiple edges. Let $\lambda(i)$ and $\lambda(\{i,j\})$ denote the color of a vertex and an edge, respectively. Let $\Gamma(\cG)$ be the group of graph automorphisms of~$\cG$ that preserve its vertex and edge colors. 
\begin{definition}\label{def:rcop}
The colored graph~$\cG$ is an \emph{RCOP graph} if 
\begin{enumerate}
\item the sets of vertex and edge colors are disjoint,
\item for any  vertices~$u,v \in [n]$ of the same color there is some~$\gamma \in \Gamma(\cG)$ with~$\gamma(u)=v$, and
\item for any edges~$e,f \in E(\cG)$ of the same color there is some~$\gamma \in \Gamma(\cG)$ with~$\gamma(e)=f$.
\end{enumerate}
\end{definition}

In the present work, we only consider \emph{connected} RCOP block graphs. Connectedness is not typically required in the graphical models literature. However, disconnected RCOP graphs have a relatively simple structure; two connected components either have no color in common or are isomorphic to each other. 

In what follows we define RCOP models and present examples and their important properties. 
RCOP models were introduced in \cite{hojsgaard2008} and have since been studied in works such as \cite{gehrmann2012,uhler2011}.

\begin{figure}[h]
    \centering
    \includegraphics[scale=0.5]{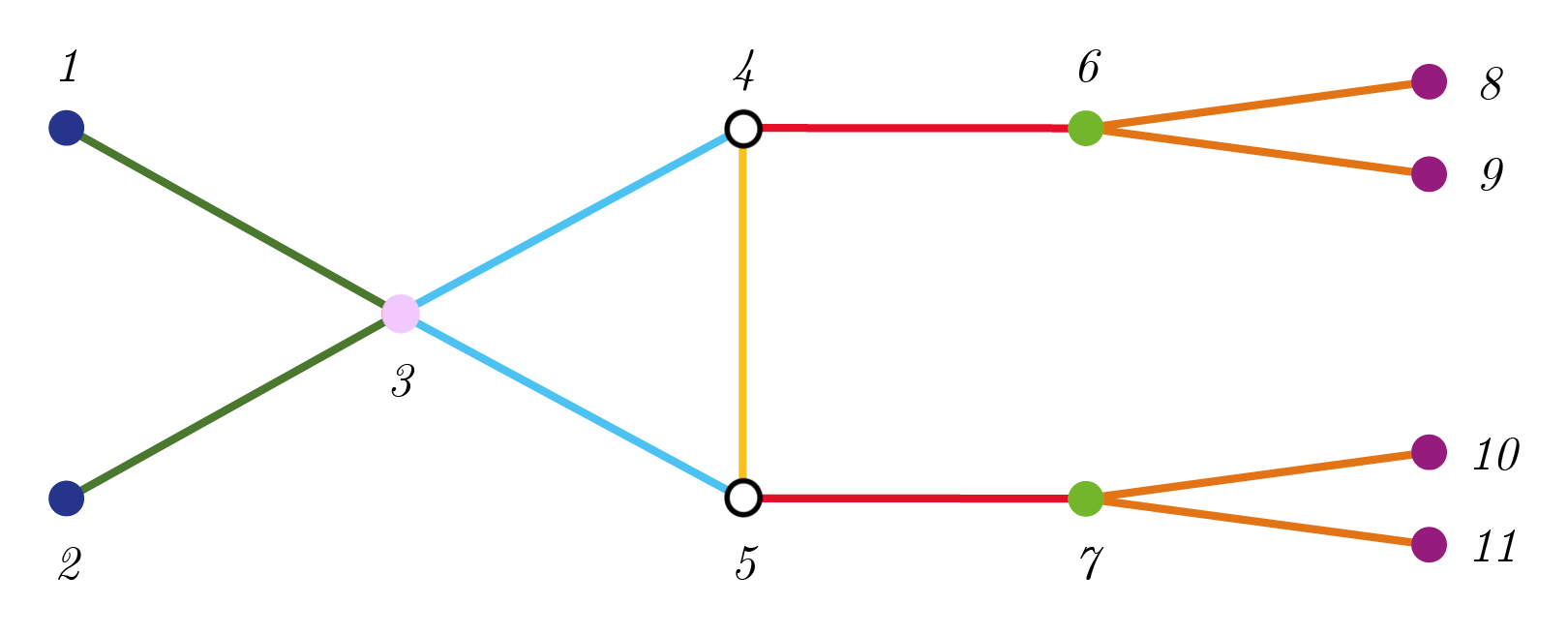}
    \caption{Example of an RCOP graph}
    \label{fig:RCOP}
\end{figure}

\begin{example}\label{ex:RCOP}
The colored graph in \Cref{fig:RCOP} is RCOP with automorphism group~\(\Gamma(\Gcal) = \langle (1 \ 2), (4\ 5)(6\ 7)(8\ 10)(9\ 11), (8\ 9), (10\ 11)\rangle.\) 
If we  change the color of one vertex in the graph, say, take~$\lambda(2)$ to be a new color, then the resulting graph is not RCOP anymore. This is because any~$\gamma$ in~$\Gamma$ with~$\gamma(\{1,3\})=\{2,3\}$ must have~$\gamma(3)=3$ and~$\gamma(1)=2$. This is not possible since vertices~$1$ and~$2$ have different colors.  
\end{example}

Let $\cG$ be a  RCOP graph on $n$ vertices.  Consider the linear space~$\mathcal{L}_{\cG}$ of symmetric matrices~$K=(k_{ij})$ in~$\mathbb{R}^{n\times n}$ that satisfy the following constraints:
\begin{enumerate}
     \item~$k_{ij}=0$ if~$\{i,j\}$ is not an edge in~$\Gcal$,
     \item~$k_{ii}=k_{jj}$ if~$\lambda(i)=\lambda(j)$ for vertices~$i,j$ in~$[n]$, and 
     \item~$k_{ij}=k_{uv}$ if~$\lambda(\{i,j\})=\lambda(\{u,v\})$ for edges~$\{i,j\}$ and~$\{u,v\}$ in~$\Gcal$.
 \end{enumerate}
 \begin{definition}\label{def:RCOPModel}
The \emph{RCOP model}~$\Mcal_{\Gcal}$ is the set of all multivariate Gaussian distributions on random variables $X_1, \dots, X_n$ with mean~$\mathbf{0}$ and positive definite concentration matrix~$K$ in~$\Lcal_{\Gcal}$. The inverse linear space, $\Lcal_{\Gcal}^{-1}=\overline{\{\Sigma \in \mathbb{R}^{n\times n} \colon \Sigma^{-1}\in \Lcal_{\Gcal}\}}$ is the Zariski closure of the set of covariance matrices for~$\Mcal_{\Gcal}$. 
\end{definition}

The entries of the positive definite covariance matrix~$\Sigma \in \Lcal_{\Gcal}^{-1}$ are the covariances between the random variables; that is,~$\sigma_{ij} = \mathrm{Cov}[X_i,X_j] = \mathbb{E}[X_iX_j] - \mathbb{E}[X_i]\mathbb{E}[X_j]$.
The entries of~$K$ are called the concentrations of the model and are useful for understanding the conditional independence constraints on the model. In particular, when~$k_{ij}=0$, this means that~$X_i$ and~$X_j$ are independent given all other random variables.

In this paper we determine conditions on an RCOP graph~$\cG$ under which~$\Lcal_{\cG}^{-1}$ is a toric variety.  We do this by proving that in these instances,~$\Lcal_{\cG}^{-1}$ is the variety of a toric ideal. Recall an ideal is \emph{toric} if and only if it is prime and generated by binomials. Equivalently, an ideal is toric if and only if it is the kernel of a monomial map from a polynomial ring to a Laurent polynomial ring. We recommend \cite{herzog2018binomial} for background on toric ideals and varieties.

The inverse linear space $\Lcal_{\cG}^{-1}$ arises as the variety of the kernel  of the rational map  \begin{align} \label{eq:minors}
\rho_\cG: \mathbb{R}[\Sigma] \to \mathbb{R}(K), \ 
\rho_\cG(\sigma_{ij})= \frac{(-1)^{i+j}K_{ij}}{\det (K)}, 
\end{align}
where~$K_{ij}$ is the~$ij$-th minor of the symmetric matrix~$K$. We denote the kernel of this map by~$I_{\cG}$ and refer to it as the  \emph{vanishing ideal of the RCOP model}~$\Mcal_\cG$. The ideal~$I_{\Gcal}$ is prime as the kernel of a rational map. We will search for conditions on $\cG$ under which $I_\cG$ is a binomial ideal.

\begin{figure}[h]
     \includegraphics[scale=0.7]{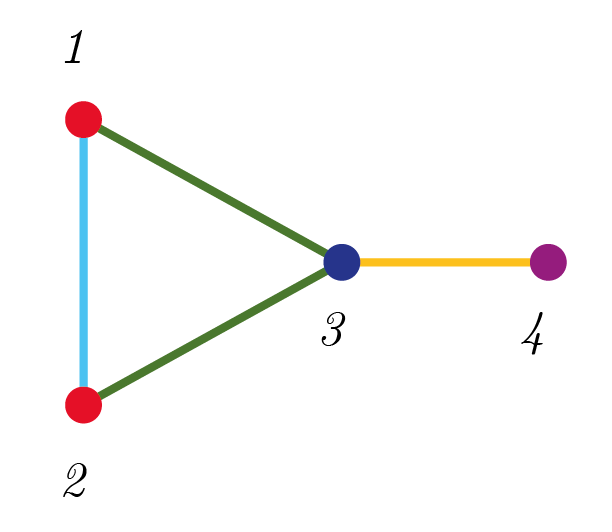}
    \hspace{1cm}
       \includegraphics[scale=0.6]{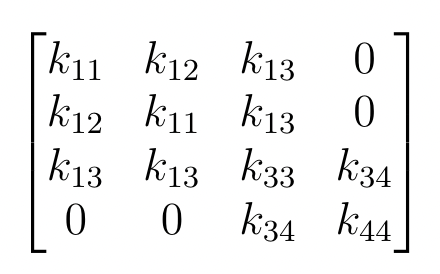}
     \hspace{1cm}
     \includegraphics[scale=0.7]{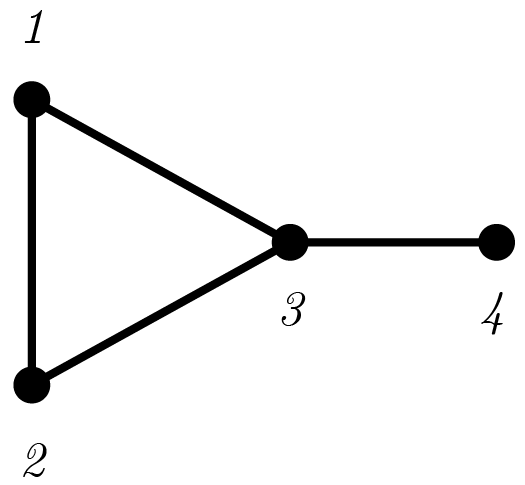}
    \caption{An RCOP graph~$\cG$, the form of its associated concentration matrix, and its underlying graph~$G$.}
    \label{fig:main}
\end{figure}

\begin{example}\label{ex:main}
(a) The colored graph $\cG$ in \Cref{fig:main} is RCOP. 
We used Macaulay2 \cite{M2} to compute the kernel of the map $\rho_{\Gcal}$ in \Cref{eq:minors} and obtained that~$I_{\cG}=\langle \sigma_{11}-\sigma_{22},\sigma_{13}-\sigma_{23},\sigma_{14}-\sigma_{24},\sigma_{24}\sigma_{33}-\sigma_{23}\sigma_{34}\rangle$.
Hence~$\Lcal^{-1}_{\cG}$ is the set of symmetric matrices~$\Sigma=(\sigma_{ij})\in \mathbb{R}^{n\times n}$ that vanish on~$I_\cG$. We will later see that this ideal is the sum of  the toric ideals~$I_{G}=\langle \sigma_{13}\sigma_{34}-\sigma_{14}\sigma_{33} ,\sigma_{23}\sigma_{34}-\sigma_{24}\sigma_{33} ,\sigma_{14}\sigma_{23}-\sigma_{13}\sigma_{24} \rangle$ of the uncolored graph~$G$, and  ~$I_{\Bar{\cG}}=\langle\sigma_{11}-\sigma_{22},\sigma_{13}-\sigma_{23},\sigma_{14}-\sigma_{24} \rangle$ of the completion graph~$\Bar{\cG}$ of~$\cG$, introduced in \Cref{def:completion} and visualized in \Cref{fig:completion}.

(b) 
\cite[Example~5.1 (\emph{Frets’ heads})]{SU10}\label{ex:Frets_heads}
Let $G$ be the $4$-cycle with edges
~$\{1,2\}, \{2,3\}, \{3,4\}$, and $\{1,4\}$.
Let $\Gcal$ have coloring given by ~$\lambda(1)=\lambda(2)$,~$\lambda(3)=\lambda(4)$, and~$\lambda(\{1,4\})=\lambda(\{2,3\})$. The edges $\{1,2\}$ and $\{3,4\}$ have distinct colors. This is an RCOP graph with~$\Gamma(\cG)=\langle (1\ 2)(3 \ 4)\rangle$.  Its vanishing ideal~$I_{\cG}=\langle \sigma_{11}-\sigma_{22}, \sigma_{33}-\sigma_{44}, \sigma_{23}-\sigma_{14}, \sigma_{13}-\sigma_{24}, \sigma_{23}\sigma_{24}-\sigma^3_{24}-\sigma_{22}\sigma_{23}\sigma_{34}+\sigma_{12}\sigma_{24}\sigma_{34}- \sigma_{12}\sigma_{23}\sigma_{44}+\sigma_{22}\sigma_{24}\sigma_{44}\rangle$ is not toric. 
\end{example}

We say that the colored graph~$\Gcal$ is \emph{edge regular} if for any two edges~$\{u,v\}$ and~$\{u',v'\}$ with~$\lambda(\{u,v\}) = \lambda(\{u',v'\})$, we have that~$\{\lambda(u), \lambda(v)\} = \{ \lambda(u'), \lambda(v')\}$. In other words, edges of the same color connect vertices of the same color. We say that~$\Gcal$ is \emph{vertex regular} if for each edge color~$\epsilon$ and pair of vertices~$u, v$ with~$\lambda(u) = \lambda(v)$, the vertices~$u$ and~$v$ are adjacent to the same number of edges of color~$\epsilon$.

\begin{lemma}\cite[Prop. 2 \& Prop. 8]{gehrmann2012}\label{lemma:regularityRCOP} 
If~$\Gcal$ is an RCOP graph, then~$\Gcal$ is edge regular and vertex regular.
\end{lemma}

One can also distinguish if a colored graph is not RCOP by comparing the neighborhoods of vertices of the same color in the graph. For a vertex~$c$ in a graph $G$, let~$V(c)$  be the set containing $c$ and all vertices adjacent to~$c$. The \emph{neighborhood} of~$c$, denoted $\Ne(c)$, is the subgraph of $G$ induced by $V(c)$.

\begin{prop}
\label{prop:RCOPneighb}
Let $c$ and $d$ be two vertices of the same color in the RCOP graph~$\cG$. Then any $\gamma\in \Gamma(\cG)$ such that $\gamma(c)=d$ has $\gamma(\Ne(c))=\Ne(d)$. In particular,~$c$ and~$d$ have isomorphic neighborhoods. 
\end{prop}

\begin{proof}
Since~$\lambda(c)=\lambda(d)$ there is at least one automorphism~$\gamma\in\Gamma(\cG)$ such that~$\gamma(c)=d$. Take any vertex~$v\in V(c)$. Since~$v$ is adjacent to~$c$, its image~$\gamma(v)$ must be adjacent to~$d$. Hence,~$\gamma(v)$ is in~$V(d)$. Similarly, take~$v\in V(d)$. Since~$v$ is adjacent to~$d$, its image under~$\gamma^{-1}\in \Gamma(\cG)$ is adjacent to~$c$. Hence,~$\gamma^{-1}(v)\in V(c)$. As result~$\gamma(V(c))=V(d)$. Since~$\gamma$ is a graph automorphism we must have~$\gamma(\Ne(c))=\Ne(d)$, as desired. 
\end{proof}

\subsection{Markov Bases}\label{subsec:Markov}
In the present work, we describe a family of RCOP models whose vanishing ideals are toric and give a generating set, or \emph{Markov basis}, for the vanishing ideal in these cases. In this section, we review the relevant theory of Markov bases. Thorough introductions to Markov basis theory can be found in \cite{markovbases2012,herzog2018binomial}. 

Let~$A \in \Z^{d \times r}$ be an integer matrix.
For each integer vector~$b$ in~$\mathbb{Z}^{r}$,
we write~$b = b^+ - b^-$ uniquely, where~$b^+$ and~$b^-$ are non-negative integer vectors with disjoint support.
Let~$\mathcal{B}=\{b_1,\dots,b_s\}$ be vectors in the integer kernel of~$A$.
Then~$\mathcal{B}$ is a \emph{Markov basis} for~$A$
if for all~$u \in \ker_{\Z}(A)$,
there is a sequence~$b_{i_1},\dots,b_{i_m}$ of vectors in~$\mathcal{B}$ such that
\begin{itemize}
    \item~$u^+ + b_{i_1} + \dots + b_{i_j}$ is non-negative for all~$j \leq m$, and 
    \item~$u^+ + b_{i_1} + \dots + b_{i_m} = u^-$.
\end{itemize}
The vectors in~$\mathcal{B}$ are called Markov moves. 

 The matrix~$A$ defines the toric ideal~$I(A)$ as the kernel of the monomial map \[\mathbb{R}[x_1,\dots,x_r]\rightarrow \mathbb{R}[y^{\pm}_1,\dots,y^{\pm}_r],\ x_i\xrightarrow[]{} y_1^{a_{1,i}}\dots y_r^{a_{r,i}}. \]
 
The following theorem, known as the Fundamental Theorem of Markov Bases, describes the relationship between Markov bases for~$A$ and the toric ideal~$I(A)$. For a vector~$b\in \mathbb{Z}_{\geq 0}^r$, we will use~$\mathbf{x}^b$ to denote the monomial~$x_1^b\cdots x_r^{b_r}$. 
\begin{theorem}\cite[Theorem 3.1]{diaconis1998algebraic}\label{thm:markov}
A set of integer vectors~$\mathcal{B} \subset \ker_{\Z}(A)$ is a Markov basis for~$A$ if and only if
$\{ \bfx^{b^+} - \bfx^{b^-} \mid b \in \mathcal{B} \}$ is a generating set for the toric ideal~$I(A)$.
\end{theorem}

This theorem establishes a bijection between Markov bases for~$A$ and generating sets of the toric ideal~$I(A)$. Thus we often refer to a generating set of~$I(A)$ as a Markov basis for the ideal and a binomial in $I(A)$ as a Markov move.

\section{Block Graphs and the Shortest Path Map}\label{sec:BlockGraphs}
The toric vanishing ideals of uncolored graphical models on \emph{block graphs} were studied extensively in \cite{misra2019gaussian}. 
In this section, we define block graphs and outline the known results related to their structure and vanishing ideals.

Let~$\cG$ be a connected RCOP graph on $n$ vertices and let~$G$ be its underlying uncolored graph. 
Let~$A,B,$ and~$C$ be disjoint subsets of~$[n]$ such that~$[n]=A\cup B\cup C$.  We say that~$C$ \emph{separates}~$A$ and~$B$ 
if for any two vertices~$a\in A$ and~$b\in B$, all paths from~$a$ to~$b$ contain a vertex in~$C$. 
\begin{definition} \label{def:block_graph}
The graph~$G$ is a \emph{block graph} (or a \emph{$1$-clique sum of complete graphs}) if there exists a partition~$(A,B,\{c\})$ of~$[n]$  for some vertex~$c$ in~$[n]$ such that  the set~$\{c\}$ separates~$A$ and~$B$, and the subgraphs induced by~$A\cup \{c\}$ and~$B\cup \{c\}$ are either complete  or block graphs.
\end{definition}

\Cref{fig:RCOP} and \Cref{fig:main} are examples of RCOP block graphs. In \cite{misra2019gaussian}, the authors show that the vanishing ideal of the set of covariance matrices arising from an uncolored Gaussian graphical model on a block graph is toric. Before we present this theorem, we first outline some useful facts about block graphs. The following proposition follows directly from the definition of a block graph.

\begin{prop}\label{prop:NoCycles}
Let $G$ be a block graph and~$D$ a cycle in $G$. Then the subgraph of $G$ induced by $D$ is a complete graph.  In particular, block graphs are \emph{chordal}. 
\end{prop}

Next we restate the useful result that there exists a unique shortest path between any two vertices in a block graph.

\begin{prop}\cite[Proposition 2]{misra2019gaussian}\label{prop:UniqueShortestPath}
Let~$G$ be a connected block graph. Then for any two vertices~$u$,~$v$ in the graph, there is 
a unique shortest path in~$G$, denoted~$u\lra v$, that connects~$u$ to~$v$. 
\end{prop}

Let $\lambda(\Gcal)$ denote the set of all colors of edges and vertices of the RCOP block graph $\Gcal$. For each pair of vertices~$u,v$ in~$\Gcal$, let~$E(u\lra v)\coloneqq \{e \mid e\in u\lra v\}$ be the set of edges used by this path, and let~$\Lambda(u\lra v)$ denote the multiset containing the colors of vertices~$u$ and~$v$, and the edge colors on the unique shortest path from~$u$ to~$v$; that is,~$\Lambda(u \leftrightarrow v) \coloneqq \{\{\lambda(u), \lambda(v)\}\} \cup \{\{\lambda(e) \mid e \in u\leftrightarrow v\}\}$ as a multiset. 

Assign to each color~$\ell$ in the  set~$\lambda(\cG)$ the parameter~$t_{\ell}$. By \Cref{prop:UniqueShortestPath}, we can define a map that sends each variable~$\sigma_{ij}$ of~$\mathbb{R}[\Sigma]$ to a product of parameters determined by the colors along the shortest path~$i\lra j$.

\begin{definition}\label{def:shortestPathMapColored}
The \emph{shortest path map} of an RCOP block graph~$\cG$ is the monomial map 
\begin{eqnarray*}
    \varphi_{\cG}  :  \mathbb{R}[\sigma_{ij} \mid i,j \in [n], i \leq j] &\rightarrow& \mathbb{R}[ t_{\ell} \mid   \ell\in \lambda(\cG)], \\
    \sigma_{ij} & \mapsto & t_{\lambda(i)} t_{\lambda(j)} \prod_{e \in E(i\leftrightarrow j)} t_{\lambda(e)}.
\end{eqnarray*}

\end{definition}
We refer to the exponent matrix defining~$ \varphi_{\cG}$ as~$A_\cG$. The kernel of~$\varphi_{\cG}$ is the toric ideal with binomials of form 
$ \sigma_{i_1j_1}\sigma_{i_2j_2}\dots \sigma_{i_dj_d}-\sigma_{k_1l_1}\sigma_{k_2l_2}\dots \sigma_{k_dl_d}$
where the multisets~$L=\{\{i_1\lra j_1,\dots, i_d\lra j_d\}\}$ and~$R=\{\{k_1\lra l_1,\dots, k_d\lra l_d\}\}$ have the same multiset of edge colors and endpoints, i.e.~$\Lambda(L)=\Lambda(R)$. 
Their associated Markov moves are integer vectors~$(b_{uv})_{1\leq u\leq v\leq n}\in \mathbb{Z}^{\binom{n}{2}+n}$ with entry~$b_{uv}$ equal to the number of times that the shortest path~$u\lra v$ appears in~$L$ if~$b_{uv}\geq 0$, and equal to the number of times that~$u\lra v$ appears in~$R$ if~$b_{uv}\leq 0$. To simplify notation, we denote the Markov move associated such a binomial by~$(i_1,j_1)(i_2,j_2)\dots (i_d,j_d) -(k_1,l_1)(k_2,l_2)\dots (k_d,l_d)$. 

\begin{example}\label{ex:ShortestPath}
The exponent matrix for  the shortest path map of the RCOP block graph~$\cG$ in \Cref{fig:main}  is
\[
A_{\cG}=
\begin{blockarray} {*{11}{c}}
& 11 & 12 & 13 & 14 & 22 & 23 & 24 & 33 & 34 & 44 \\
\begin{block}{c[*{10}{c}]}
{\red{\boldsymbol{\lambda(1)}}} & 2 & 2 & 1 & 1 & 2 & 1 & 1 & 0 & 0 & 0 \\
{\blue{\boldsymbol{\lambda(3)}}} & 0 & 0 & 1 & 0 & 0 & 1 & 0 & 2 & 1 & 0 \\
{\color{violet}{\boldsymbol{\lambda(4)}}} & 0 & 0 & 0 & 1 & 0 & 0 & 1 & 0 & 1 & 2 \\
{\color{cyan}{\boldsymbol{\lambda(\{1,2\})}}} & 0 & 1 & 0 & 0 & 0 & 0 & 0 & 0 & 0 & 0 \\
{\color{ForestGreen}{\boldsymbol{\lambda(\{1,3\})}}} & 0 & 0 & 1 & 1 & 0 & 1 & 1 & 0 & 0 & 0 \\
{\color{Dandelion}{\boldsymbol{\lambda(\{3,4\})}}} & 0 & 0 & 0 & 1 & 0 & 0 & 1 & 0 & 1 & 0 \\
\end{block}
\end{blockarray} \hspace{.2cm}.
\]
The rows of this matrix are indexed by the colors of the vertices and edges. The kernel of~$\varphi_{\cG}$ has a Markov basis~$\mathcal{B}_{\cG}=\{(1,1)-(2,2), (1,3)-(2,3),(1,4)-(2,4), (2,3)(3,4)-(2,4)(3,3)\}$ which coincides with the generating set of the ideal~$I_{\cG}$ computed in \Cref{ex:main}. 
\end{example}

The uncolored version of this map, wherein each vertex and edge has a different parameter, was previously introduced by Misra and Sullivant in \cite[Definition 4]{misra2019gaussian}. In the case of an uncolored graph~$G$, a move of the form~$L-R$ is in the ideal~$I_G$ if and only if the multiset of edges and endpoints used in all paths in~$L$ is equal to that of all paths in~$R$. Now we are ready to state the main known result on standard Gaussian graphical models on block graphs.  
\begin{theorem}[Theorem 1 \& Theorem 5, \cite{misra2019gaussian}]\label{thm:UncoloredMarkovBasis}
Let~$G$ be a block graph.
Then~$I_G$ is the toric ideal~$\ker(\varphi_G)$ and it is  generated by the quadratic binomials 
\begin{align}\label{eq:gensG}
   \{\sigma_{ij} \sigma_{k\ell} - \sigma_{ik}\sigma_{j\ell} \mid E(i \leftrightarrow j) \cup E(k \leftrightarrow \ell) = E(i \leftrightarrow k) \cup E(j \leftrightarrow \ell) \}, 
\end{align}
where each union is taken as a multiset. 
\end{theorem}

\begin{example}\label{ex:UncoloredShortestPath}

The Markov basis for the kernel of the shortest path map of the graph~$G$ in \Cref{fig:main} 
is~$\mathcal{B}_G=\{(1,3)(3,4)-(1,4)(3,3), (2,3)(3,4)-(2,4)(3,3), (1,4)(2,3)-(1,3)(2,4)\}$. This coincides with the generators of the ideal~$I_G$ computed in \Cref{ex:main}. 
\end{example}

\section{Properties of RCOP Block Graphs}\label{sec:completionsAreRCOP}
In this section, we analyse the combinatorial properties of RCOP block graphs. These results will allow us to prove that the kernel of the shortest path map from \Cref{def:shortestPathMapColored} is  the vanishing ideal of the RCOP model.

\subsection{Shortest Paths} In the first part of this section, we discuss some properties of shortest paths in RCOP block graphs. Throughout the rest of the paper, we let $\Gcal$ be a connected RCOP block graph with associated automorphism group $\Gamma(\Gcal)$.

\begin{definition}
Two paths~$u\lra v$ and~$u'\lra v'$ in a RCOP graph are \emph{combinatorially equivalent} if~$\Lambda(u\lra v)=\Lambda(u'\lra v')$. If the edge colors of the two paths appear in the same order, we say that the paths are \emph{isomorphic} to one another. 
\end{definition}
As we will soon see, these two notions  are equivalent in RCOP block graphs. 
 First we notice that an automorphism in~$\Gamma(\Gcal)$ maps each shortest path to a shortest path that is isomorphic (and hence, combinatorially equivalent) to it. 
 We denote by~$v_1\lra v_s$ the 
 shortest path that passes through  vertices~$v_1,v_2,\dots,v_s$ in the presented order.

\begin{lemma}\label{lm:image}
Let~$\cG$ be an RCOP block graph. 
The image of the shortest path~$v_1\lra v_s$ in~$\cG$ under any~$\gamma \in \Gamma(\cG)$ is the  shortest path~$\gamma(v_1)\lra \gamma(v_s)$ which is isomorphic to~$v_1\lra v_s$. 
\end{lemma}

\begin{proof}
Since~$\gamma$ sends edges of a color class to edges of the same color class, the image path~$\gamma(v_1\lra v_s)$ is isomorphic to~$v_1\leftrightarrow v_s$. 
Since $\gamma^{-1}$ is also an automorphism and  $v_1 \lra v_s$ is the shortest path between $v_1$ and $v_s$, $\gamma(v_1\lra v_s)$ is the shortest path between $\gamma(v_1)$ and $\gamma(v_s)$.
\end{proof}

The following lemma shows that shortest paths that have endpoints of the same color must have a symmetric color pattern. 

\begin{lemma}
\label{lm: symmetry}
Let~$u_1 \leftrightarrow u_s$ be a shortest path in the  RCOP block graph~$\Gcal$ with ~$\lambda(u_1) = \lambda(u_s)$. Then the sequences of edge colors~$\lambda(\{u_1,u_2\}), \lambda(\{u_2,u_3\}),\ldots, \lambda(\{u_{s-1},u_s\})$ and vertex colors $\lambda(u_1),\lambda(u_2),\ldots, \lambda(u_s)$ are symmetric.
\end{lemma}

\begin{proof}
The symmetry of the sequence of vertex colors follows from edge regularity and the symmetry of the edge colors; so it suffices to prove that the edge colors in such a path are symmetric. For the sake of contradiction, suppose there is a shortest path between vertices of the same color whose edge colors are not symmetric. Let $u_1 \lra u_s$ be the shortest path of maximal length in $\Gcal$ such that $\lambda(u_1) = \lambda(u_s)$ but $\lambda(\{u_1,u_2\}) \neq \lambda(\{u_{s-1},u_s\})$. Then there exists a $\gamma \in \Gamma(\Gcal)$ such that $\gamma(u_1) = u_s$. Note that $\gamma(u_s) \neq u_1$ since by \Cref{lm:image}, $\gamma(u_1),\dots,\gamma(u_s)$ is the shortest path from $\gamma(u_1)$ to $\gamma(u_s)$ and we have $\gamma(\{u_{s-1},u_s\}) \neq \{u_1,u_2\}$.
Moreover, we cannot have $\gamma(u_i) = u_j$ for any $i,j \in [s]$ as this would introduce a cycle and hence, an alternate shortest path between vertices of $u_1\lra u_s$ or $\gamma(u_1) \lra \gamma(u_s)$. 

Consider $u_1 \lra \gamma(u_s)$. Let $i$ be maximal such that $u_i$ lies on this path and let $j$ be minimal such that $\gamma(u_j)$ lies on this path. We claim that either we have $i=s, j=1$ or $i=s-1, j=2$. Suppose that $i \neq s$. Then $j \neq 1$ since $\gamma(u_1) = u_s$ does not lie on $u_1 \lra \gamma(u_s)$. Thus the concatenation of paths $u_i \lra u_s \lra \gamma(u_j) \lra u_i$ is a cycle, which is induces a complete graph by \Cref{prop:NoCycles}.
But $u_s$ is not adjacent to $u_i$ for $i < s-1$ or $\gamma(u_j)$ for $j > 2$. Thus $i = s-1, j = 2$.

Thus the shortest path $u_1 \lra \gamma(u_s)$ is either of the form $u_1 \lra u_s \lra \gamma(u_s)$ or $u_1 \lra u_{s-1} \lra \gamma(u_2) \lra \gamma(u_s)$. Both of these paths have length greater than or equal to $2s-1 >s$ and do not have symmetric edge colors. This contradicts the maximality of $u_1 \lra u_s$.
\end{proof}

Another very useful property of~$\cG$ is that no vertex color can appear more than twice in a shortest path, which we shall show in \Cref{lm:twovertices}. To prove this, we first need the following result that holds for any block graph. 

\begin{lemma} \label{lm:union}
Let~$\cG$ be a block graph. Let~$v_1\lra v_{k+1}= v_1,v_2,\dots,v_{k+1}$ and let~$v_{k}\lra v_s= v_{k},v_{k+1},\dots,v_s$ be two shortest paths that share the edge~$\{v_{k},v_{k+1}\}$. Then their union path~$v_1,v_2,\dots,v_s$ is the shortest path~$v_1\lra v_s$. 
\end{lemma}

\begin{proof}
Consider the shortest path $v_1 \lra v_s$. Let $i\leq k$ be maximal such that $v_i$ lies on $v_1 \lra v_s$. Let $j \geq k+1$ be minimal such that $v_j$ lies on $v_1 \lra v_s$. Note that $v_{k+1}$ is not adjacent to $v_a$ for any $a<k$ and $v_k$ is not adjacent to $v_b$ for any $b > k+1$. Thus by \Cref{prop:NoCycles}, the concatenation of paths $v_i \lra v_{k+1} \lra v_j \lra v_i$ does not form a cycle of length greater than or equal to three. Thus we have $i = k$ and $j = k+1$, as needed.
\end{proof}

\begin{lemma}
\label{lm:twovertices}
Any shortest path in an RCOP block graph has at most two vertices of the same color. 
\end{lemma}

\begin{proof}
In order to prove that any shortest path has at most two vertices of the same color, it is enough to show that the shortest path between any two vertices of the same color does not pass through another vertex of the same color. Let $u_1 \lra u_s$ be a path with $\lambda(u_1) = \lambda(u_s)$ that passes through another vertex of color $\lambda(u_1)$.  Let $k$ be minimal such that $1 < k < s$ and $\lambda(u_k) = \lambda(u_1)$. 
We claim that either $k=2$ or $k=3$. If $k \neq 2$, then by \Cref{lm: symmetry},  the three shortest paths $v_1\leftrightarrow v_k$, $v_k\leftrightarrow v_m$ and $v_1\leftrightarrow v_m$ are all symmetric. Thus $\lambda(v_2) = \lambda(v_{k+1})$ and $v_2 \lra v_{k+1}$ is symmetric. Hence $\lambda(v_3) = \lambda(v_k)$. Since $k$ was chosen to be minimal, we have that $k=3$, as needed. So by induction on $s$, such a path either has all vertices of the same color or alternates between two vertex colors.

\emph{Case 1:} 
Let $v_1\lra v_s$ be a path of maximal length such that all vertices are of color $\lambda(v_1)$ and all edges are of color $\lambda(\{v_1,v_2\})$.  
Take $\gamma\in \Gamma(\cG)$ such that $\gamma(\{v_{s-2},v_{s-1}\})=\{v_{s-1}, v_s\}$, and consider $\gamma(v_s)$. 

Suppose that $\gamma(v_{s-2})=v_{s-1}$ and  $\gamma(v_{s-1})=v_{s}$.  Since $\gamma$ is a graph automorphism,  $\gamma(v_s)\neq v_{s-1},v_s$ and $\gamma(v_s)$ is adjacent to $v_s$. Moreover,  $\gamma(v_s)\neq v_i$ for $1\leq i\leq s-2$, because  $v_1,\dots, v_s$ is the shortest path from $v_1$ to $v_s$. The two shortest paths $v_1\lra v_s$ and $\gamma(v_{s-2})\lra \gamma(v_s)$ satisfy the conditions of \Cref{lm:union} with common edge $\{v_{s-1},v_s\}$. Since $\lambda(\{v_{s-1},v_s\})=\lambda(\{v_s,\gamma(v_s)\})$, $v_1\lra \gamma(v_s)$ has the desired color pattern and is of length $s+1$, which contradicts the maximality of $v_1 \lra v_s$.

Suppose now that $\gamma(v_{s-2})=v_{s}$ and  $\gamma(v_{s-1})=v_{s-1}$.  The two shortest paths~$v_1\lra v_s$ and~$\gamma(v_{s-1})\lra \gamma(v_{1})$ satisfy the conditions of \Cref{lm:union} since they overlap on the edge $\{v_{s-1}, v_s\}$.
 The resulting path~$v_1\lra \gamma(v_1)$ has the desired color pattern and is of length~$2s-1$, which again contradicts the maximality of $v_1 \lra v_s$. 

\emph{Case 2:} 
Let $v_1\lra v_s$ be a shortest path of maximal length whose vertex colors  alternate between $\lambda(v_1)$ and $\lambda(v_2)$, and whose edge colors are all $\lambda(\{v_1,v_2\})$. Then $s \geq 5$.
Since~$\lambda(\{v_{s-2},v_{s-1}\})=\lambda(\{v_{s-1},v_s\})$, there is some~$\gamma\in \Gamma(\cG)$ such that~$\gamma(\{v_{s-2},v_{s-1}\})=\{v_{s-1},v_s\}$. Since~$\lambda(v_{s-2})\neq \lambda(v_{s-1})$, we must have~$\gamma(v_{s-2})=v_s$ and~$\gamma(v_{s-1})=v_{s-1}$. 
The shortest paths~$v_1\lra v_s$ and~$\gamma(v_{s-1})\lra \gamma(v_{1})$ satisfy the conditions of~\Cref{lm:union} with common edge~$\{v_{s-1},v_s\}$. Their union~$v_1\lra \gamma(v_1)$ is a shortest path of length~$2s-3>s$. 
\end{proof}

Now we are ready to prove one of the main results of this section, which implies that for RCOP block graphs, combinatorial equivalence and isomorphism of shortest paths are equivalent notions. 

\begin{prop}\label{prop:ComboEquivalent}
Any two combinatorially equivalent shortest paths in an RCOP block graph are isomorphic to each other. 
\end{prop}

\begin{proof}
Consider two combinatorially equivalent shortest paths,~$u_1\lra u_s$ and~$v_1\lra v_s$, with~$\lambda(u_1)=\lambda(v_1)$ and~$\lambda(u_s)=\lambda(v_s)$. Assume for the sake of contradiction that the two paths are not isomorphic to each other. We may also assume without loss of generality that~$\lambda(\{u_1,u_2\})\neq\lambda(\{v_1,v_2\})$ since otherwise, we could take the two paths that begin at the first vertices where their edge colors differ. Since the two paths are combinatorially equivalent, there must be an edge~$\{u_i,u_{i+1}\}$  with~$i>1$ of color~$\lambda(\{v_1,v_2\})$.
By edge regularity and \Cref{lm:twovertices}, this implies that exactly one of~$u_i$ and~$u_{i+1}$ has color~$\lambda(u_1)$.
If~$\lambda(u_{i+1})= \lambda(u_1)$,  the shortest path~$u_1\lra u_{i+1}$ is not symmetric and contradicts \Cref{lm: symmetry}. Hence,~$u_i$ must have color~$\lambda(u_1)$.  Similarly, there must be an edge~$\{v_k,v_{k+1}\}$ where~$k>1$  of color~$\lambda(\{u_1,u_2\})$  in~$v_1\lra v_s$, and  its endpoint~$v_k$ has color~$\lambda(u_1)$. By \Cref{lm: symmetry}, we have~$\lambda(\{v_{k-1},v_{k}\})=\lambda(\{v_1,v_2\})$. We consider the two cases where $k=2$ and $k >2$. 

If $k=2$, then by edge regularity, we have $\lambda(u_i) = \lambda(u_{i+1})$. But since $i \neq 1$, the vertices $u_1, u_i$ and $u_{i+1}$ are three distinct vertices of the same color on $u_1 \lra u_s$. This contradicts \Cref{lm:twovertices}.

If $k > 2$, there are two distinct edges~$\{v_1,v_2\}$ and~$\{v_{k-1},v_k\}$ of the same color in~$v_1\leftrightarrow v_s$. This implies that
there is another edge~$\{u_j,u_{j+1}\}$ of color~$\lambda(\{v_{k-1},v_{k}\})$  in~$u_1\lra u_s$ where~$j\neq 1,i$. Thus the path~$u_1\lra u_s$ contains the three distinct vertices~$u_1,u_i$ and $u_j$ of the same color, which contradicts \Cref{lm:twovertices}.
\end{proof}

We conclude this subsection with the following technical lemma which relies on the symmetry of shortest paths between vertices of the same color. This lemma will be critical to the proof of the main result in \Cref{sec:Main}.

\begin{lemma}\label{lm:i2j2existance}
Let $c \lra j$ and $c \lra l$ be two paths that intersect only at $c$. Let $a$ and $c+1$ denote the vertices that follow $c$ on $c \lra j$ and $c \lra l$, respectively. Suppose that $\lambda(\{c,a\}) \neq \lambda(\{c,c+1\})$.
Then either there is no edge in $a \lra j$ of color $\lambda(\{c,c+1\})$ or there is no edge in $c+1 \lra l$ of color $\lambda(\{c,a\})$.
\end{lemma}

\begin{figure}[ht]
    \centering
    \includegraphics[scale=.6]{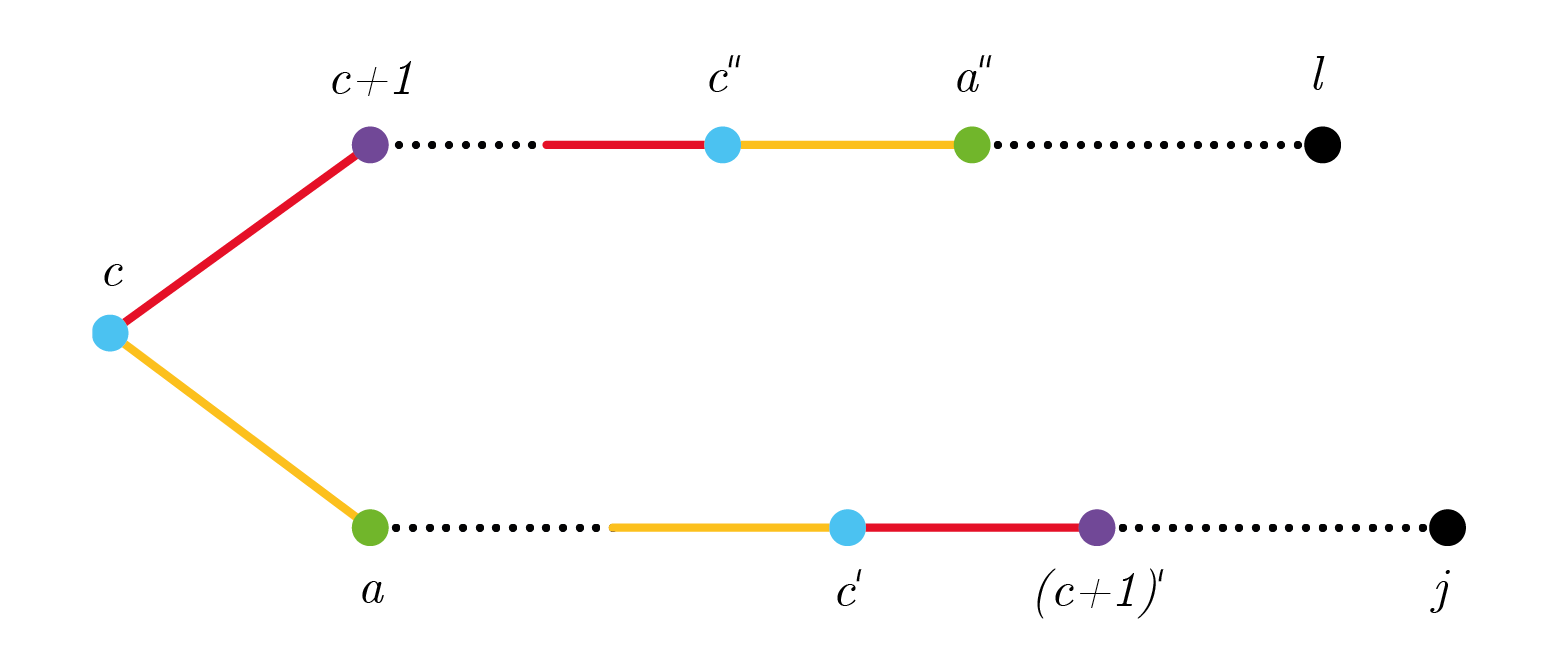}
 \caption{Paths~$c \lra j$ and~$c \lra l$. Vertices and edges marked in black have an unspecified color, whereas we use other colors to indicate that certain vertices or edges have the same image under coloring~$\lambda$. Note that the proof does not require $a$ and $c+1$ to be non-adjacent.}
 \label{fig:OriginalPaths}
\end{figure}

\begin{proof}
We must show that there cannot simultaneously be an edge in~$a\lra j$ of color~$\lambda(\{c,c+1\})$ and an edge in~$(c+1)\lra l$  of color~$\lambda(\{c,a\})$. \Cref{fig:OriginalPaths} depicts the relevant paths in~$\Gcal$.
Indeed, for the sake of contradiction, let~$\{c',(c+1)'\}$ be an edge of color~$\lambda(\{c,c+1\})$   in~$a\leftrightarrow j$ and~$\{c'',a''\}$ be of color~$\lambda(\{a,c\})$ in~$c+1\leftrightarrow l$. 
Specifically, we have~$\lambda(c') = \lambda(c'') = \lambda(c)$,~$\lambda((c+1)') = \lambda(c+1)$, and~$\lambda(a'') = \lambda(a)$. By \Cref{lm:twovertices}, $c'$ is the only other vertex of color $\lambda(c)$ on $c \lra j$, and similarly $c''$ is the only other vertex of color $\lambda(c)$ on $c \lra l$. 
Thus by \Cref{lm: symmetry} vertices~$c'$ and~$c''$ must appear before~$(c+1)'$ and~$a''$ in the paths~$(c+1)\lra (c+1)'$ and~$a\lra a''$, respectively.  The color of the other edge  adjacent to~$c''$ in~$c\leftrightarrow l$  must be~$\lambda(\{c,c+1\})$ due to the symmetry of~$c\leftrightarrow c''$. Similarly the other edge adjacent to~$c'$ in the path~$c\leftrightarrow j$ must be of color~$\lambda(\{c,a\})$. Note that we cannot have $c'=a$ or $c'' = c+1$. Indeed, if $c'=a$, this implies that $\lambda(a'') = \lambda(c'') = \lambda(c)$, which contradicts \Cref{lm:twovertices}. Thus, $\{c',c+1\}$ and $\{c'',a\}$ are not edges in $\Gcal$. This, along with the chordality of $\Gcal$, implies that the shortest path $c' \lra c''$ begins with an edge of color $\lambda(\{c,a\})$ and ends with an edge of color $\lambda(\{c,c+1\})$, which contradicts \Cref{lm: symmetry}.
\end{proof}

\subsection{Automorphisms} 

In the next chain of lemmas we explore properties of the automorphism group of an RCOP block graph. 
The results are critical to the proofs of 
\Cref{thm:rcop_complete} and \Cref{cor:CompletionGenerators}, which allow us to characterize the linear forms that vanish on~$\mathcal{L}^{-1}_{\cG}$.
First we need some notation.

\begin{definition}
 Let $c$ be a vertex of the block graph $\cG$. Consider the connected components of $\cG$ minus $c$. The $c$-\emph{components} are the subgraphs obtained by adding $c$ back into each of these connected components.
\end{definition}

For illustration, the~$3$-\emph{components} of the graph in  \Cref{fig:RCOP} are the induced subgraphs on vertex sets~$\{1,3\}$,~$\{2,3\}$ and~$\{3,4,\dots,11\}$. 
We note that if $v\neq c$, then the shortest path~$c\lra v$ is fully contained in the unique~$c$-component of~$\cG$ containing vertex~$v$.
The first lemma is an argument that we use extensively throughout the rest of this paper.

\begin{lemma}\label{lm:iso_cliques}
Let~$\{u_1,v_1\}$ and~$\{u_2,v_2\}$ be two edges of the same color contained in maximal cliques~$\cC_1$ and~$\cC_2$, respectively, in an RCOP block graph~$\cG$. Then for each~$\gamma \in \Gamma(\cG)$, if~$\gamma(\{u_1,v_1\}) = \{u_2,v_2\}$, then $\gamma(\cC_1) = \cC_2$.
In particular, any two cliques that share an edge color in an RCOP block graph are isomorphic. 
\end{lemma}
 
\begin{proof}
{As $\gamma(\{u_1,v_1\})=\{u_2,v_2\}$, $\gamma(u_1)$ can be either $u_2$ or $v_2$.} Without loss of generality, let~$\gamma(u_1)=${$u_2$} and~{$\gamma(v_1)=v_2$}. Any vertex~$w$ in~$\cC_1$ is adjacent to both~$u_1$ and~$v_1$. Its image~$\gamma(w)$ must be adjacent to both~$\gamma(u_1)$ and~{$\gamma(v_1)$}.  Since~$\cG$ is a block graph, the vertices of~$\cC_2$  are the only 
{vertices in $\cG$ that are adjacent to both $\gamma(u_1)$ and $\gamma(v_1)$.}
Hence,~$\gamma(\cC_1)\subseteq \cC_2$. For the other inclusion,  using~$\gamma^{-1}$ similarly shows that~$\gamma^{-1}(\cC_2)\subseteq \cC_1$, and so~$\cC_2\subseteq \gamma(\cC_1)$. 
\end{proof}

Our next aim is to prove \Cref{cor:completionRCOP}, which connects any two isomorphic shortest paths by an automorphism in $\Gamma(\cG)$. To do so, we must understand the behavior of $\Gamma(\cG)$ on the components of vertices of the same color.

\begin{lemma}\label{lm:alphabeta}
Let $\cG$ be an RCOP block graph. 
\begin{enumerate}[label=(\arabic*),itemsep=5pt]
    \item Let~$\{u,v\}$ be an edge in~$\cG$ with~$\lambda(u)=\lambda(v)$. Take~$\cG_u$ to be the union of all of the~$u$-components except for the one containing~$v$. Similarly, let~$\cG_v$ be the union of all of the~$v$-components except for the one containing~$u$. There exists an automorphism~$\alpha$ in~$\Gamma(\cG)$ such that~$\alpha(\cG_u)=\cG_v$,~$\alpha(\cG_v)=\cG_u$ and~$\alpha(w)=w$ for all vertices~$w$ not in~$\cG_u\cup \cG_v$. In particular,~$\alpha(u)=v$ and~$\alpha(v)=u$. \label{lm:alpha} 
    
    \item 
Let~$\{c,u\}$ and~$\{c,v\}$ be two edges of the same color in~$\cG$. Let~$H_u$ and~$H_v$ be the~$c$-components {containing}~$u$ and~$v$, respectively. There exists an automorphism~$\beta$ in~$\Gamma(\cG)$ such that~$\beta(H_u)=H_v$,~$\beta(H_v)=H_u$, and~$\beta(w)=w$ for all vertices~$w$ not in~$H_u\cup H_v$. In particular,~$\beta(u)=v$,~$\beta(v)=u$, and~$\beta(c)=c$.\label{lm:beta}
\end{enumerate}
\end{lemma}

\begin{proof}
(1) Since the two vertices~$u$ and~$v$ share the same color, there is some~$\gamma$ in~$\Gamma(\cG)$ such that~$\gamma(u)=v$. Take a vertex~$w$ in~$\cG_u$. Then,~$w\lra u$ is a shortest path fully contained in one of the~$u$-components in~$\cG_u$. Let~$u_1 \in \cG_u$ be the vertex in~$w\lra u$ adjacent to~$u$. Suppose for contradiction that~$\gamma(w)$ is in the~$v$-component containing~$u$.  Then~$\gamma(u_1)$ is one of the vertices in the unique maximal clique~$\cC$ containing the edge~$\{u,v\}$. By \Cref{lm:iso_cliques}, the maximal clique~$\cC_1 \subset \cG_u$ containing~$\{u_1,u\}$ is isomorphic to~$\cC$ via~$\gamma$. In particular, since~$\cC$ contains at least two vertices of color~$\lambda(u)$,  there is a vertex~$u_2\neq u$ in~$\cC_1$ of color~$\lambda(u)$. The shortest path $u_2 \lra v$  passes through~$u_2$,~$u$ and~$v$ and thus contains three vertices of the same color, which contradicts \Cref{lm:twovertices}. Hence,~$\gamma(\cG_u)$ must be a subset of~$\cG_v$.
Using~$\gamma^{-1}$ we see that~$\gamma^{-1}(\cG_v)\subseteq \cG_u$, and hence~$\cG_v\subseteq \gamma(\cG_u)$. Hence~$\gamma(\cG_u)=\cG_v$. 
Lastly, use~$\gamma$ to construct the automorphism
\begin{align*}
\alpha(i)=\begin{cases}
\gamma(i) & \text{ for } i\text{ a vertex in } \cG_u,\\
\gamma^{-1}(i) & \text{ for } i\text{ a vertex in } \cG_v, \text{ and }\\
i & \text{ otherwise}.
\end{cases}
\end{align*}
This is indeed an automorphism as the only vertices that do not lie in $\cG_u$ or $\cG_v$ all belong to the clique $\cC$.
 The vertex and edge colors in~$\cG$ are invariant under~$\alpha$ by definition, which makes~$\alpha$ an element of~$\Gamma(\cG)$. Thus, the automorphism~$\alpha$ satisfies all the desired properties of the lemma.\\
 
 (2) We will first show that there exists some automorphism~$\gamma$ in~$\Gamma(\cG)$ such that~$\gamma(u)=v$ and~$\gamma(c)=c$. Indeed, since edges~$\{c,u\}$ and~$\{c,v\}$ have the same color, there is some~$\delta$ in~$\Gamma(\cG)$ such that~$\delta(\{c,u\})=\{c,v\}$. If $\delta(u)=v$, then take $\gamma\coloneqq \delta$. Otherwise, $\delta(u)=c$ and~$\delta(c)=v$. Then, $\lambda(c)=\lambda(u)=\lambda(v)$ since $\cG$ is RCOP. By \cref{lm:twovertices}, $u$ and $v$ must be adjacent to each other. Take~$\alpha$ as in \Cref{lm:alphabeta}\ref{lm:alpha} applied to vertices~$c$ and~$v$. The composition~$\gamma\coloneqq \alpha \circ \delta$ has~$\gamma (c)=c$ and~$\gamma(u)=v$, as desired.

We will construct the desired map~$\beta$ from~$\gamma$ in the previous paragraph. Before this, we must show that~$\gamma(H_u)=H_v$. Take a vertex~$w \in H_u$. The shortest path~$c\lra w$ lies completely in~$H_u$. Consider the first edge $\{c,u_1\}$ that belongs to $c\lra w$.  Since this edge is in the clique containing $\{c,u\}$, we must have $u_1=u$ or $u_1$ is adjacent to $u$. Either way, $c,u$ and $u_1$ are in the same maximal clique. By \Cref{lm:iso_cliques}, the first edge $\gamma(\{c,u_1\})=\{c,\gamma(u_1)\}$ belonging to $c\lra \gamma(w)$ must be in the unique clique containing $\{c,v\}$. Hence, the entire path $c\lra \gamma(w)$ is in $H_v$. 
Similarly we see that~$\gamma^{-1}(H_v)\subseteq H_u$, and so~$H_v\subseteq \gamma(H_u)$.

Lastly,  define the automorphism
\begin{align}
\label{eq:alpha}
\beta(i)=\begin{cases}
\gamma(i) & \text{ for } i\text{ a vertex in } H_u,\\
\gamma^{-1}(i) & \text{ for } i\text{ a vertex in } H_v, \text{ and }\\
i & \text{ otherwise}.
\end{cases}
\end{align}
By the construction of~$\beta$, the vertex and edge colors in~$\cG$ are invariant under~$\beta$. So, the automorphism $\beta \in \Gamma(\cG)$ satisfies the desired properties. 
\end{proof}

We advise the reader to compare the generators of $\Gamma(\cG)$ in  \Cref{ex:RCOP} with the automorphisms in  (\ref{eq:alpha}). 
We end this section by proving that for any  two combinatorially equivalent (and hence, isomorphic) paths in an RCOP block graph, there is an automorphism in~$\Gamma(\cG)$ that maps one path to the other.

\begin{prop}\label{cor:completionRCOP}
Let~$u_1\lra u_s$ and~$v_1\lra v_s$ be two isomorphic shortest paths in the RCOP block graph~$\cG$ with  $\lambda(u_i)=\lambda(v_i)$ for all~$i$ in~$[s]$. Then there is an automorphism in~$\Gamma(\cG)$ that maps each $u_i$ to $v_i$.
\end{prop}

\begin{proof}
We use induction on the number of vertices in the two paths. In the base case when~$s=2$, these paths are just the two edges~$\{u_1,u_2\}$ and~$\{v_1,v_2\}$  of the same color in~$\cG$. 
By definition of an RCOP graph, there is some~$\delta$ in~$\Gamma(\cG)$  with~$\delta(\{u_1,u_2\})=\{v_1,v_2\}$.  If~$\delta(u_1)=v_1$ and~$\delta(u_2)=v_2$, then $\delta$ is the desired automorphism. Otherwise, given that~$\delta$ preserves the color of vertices, we have~$\lambda(v_1)=\lambda(v_2)$. By \Cref{lm:alphabeta}\ref{lm:alpha}, there exists an automorphism~$\alpha$ in~$\Gamma(\cG)$ such that~$\alpha(v_1)=v_2$ and ~$\alpha(v_2) = v_1$. The composition~$\alpha\circ \gamma$ has the desired property. 
 
 Assume now that for any two combinatorially equivalent shortest paths that use~$s$ vertices, there is an automorphism in~$\Gamma(\cG)$ mapping one path to the other. Let $u_1 \lra u_{s+1}$ and $v_1 \lra v_{s+1}$ be combinatorially equivalent shortest paths with $s+1$ vertices. By \Cref{prop:ComboEquivalent}, these paths are isomorphic; so the shortest paths~$u_1\lra u_{s}$ and~$v_1\lra v_{s}$ are also isomorphic. By the inductive hypothesis, there is some~$\delta$ in~$ \Gamma(\cG)$ such that~$\delta(u_i)=v_i$ for~$i=1,\dots,s$.  
The edges~$\{v_{s},v_{s+1}\}$ and~$\{v_{s}, \delta(u_{s+1})\}$ satisfy conditions of  \Cref{lm:alphabeta}\ref{lm:beta}. Hence, there is an automorphism~$\beta$ in~$\Gamma(\cG)$ with~$\beta(v_i)=v_i$ for~$i=1,\dots,s$ and~$\beta(\delta(u_{s+1}))=v_{s+1}$. Hence we have  $\beta \circ \delta(u_i)=v_i$ for~$i=1,\dots,s+1$, as desired. 
\end{proof}

\section{Completions of RCOP Block Graphs and Jordan Algebras}\label{sec:Jordan}

In this section we investigate the linear forms in the vanishing ideals of RCOP models. To achieve this, we define completions of RCOP block graphs, as mentioned in \Cref{thm:IntroMain}, and connect them to Jordan algebras.
 
\begin{definition}\label{def:completion}
Let~$\Gcal$ be an RCOP block graph on vertex set~$[n]$ with edge set~$E$. The \emph{completion} of~$\Gcal$, denoted~$\Gbar$, is the complete graph on~$[n]$ with the coloring~$\overline{\lambda}$ defined as follows:
\begin{enumerate}
    \item~$\overline{\lambda}(v) = \lambda(v)$ and~$\overline{\lambda}(e) = \lambda(e)$ for all~$v \in [n]$ and~$e \in E$ and 
    \item the edges in~$\Gbar$ that are not in~$E$ have \emph{new} colors assigned to them such that  any two edges~$\{u,v\}, \{u',v'\}$ of~$\Gbar$ that are not in~$E$ share the same color if and only if~$\Lambda({u\lra v}) = \Lambda({u'\lra v'})$.
\end{enumerate}
\end{definition}

\begin{figure}[h]
    \centering
    \includegraphics[scale=0.7]{Figures/figure1_rcop.png}
    \hspace{1cm}
    \includegraphics[scale=0.3]{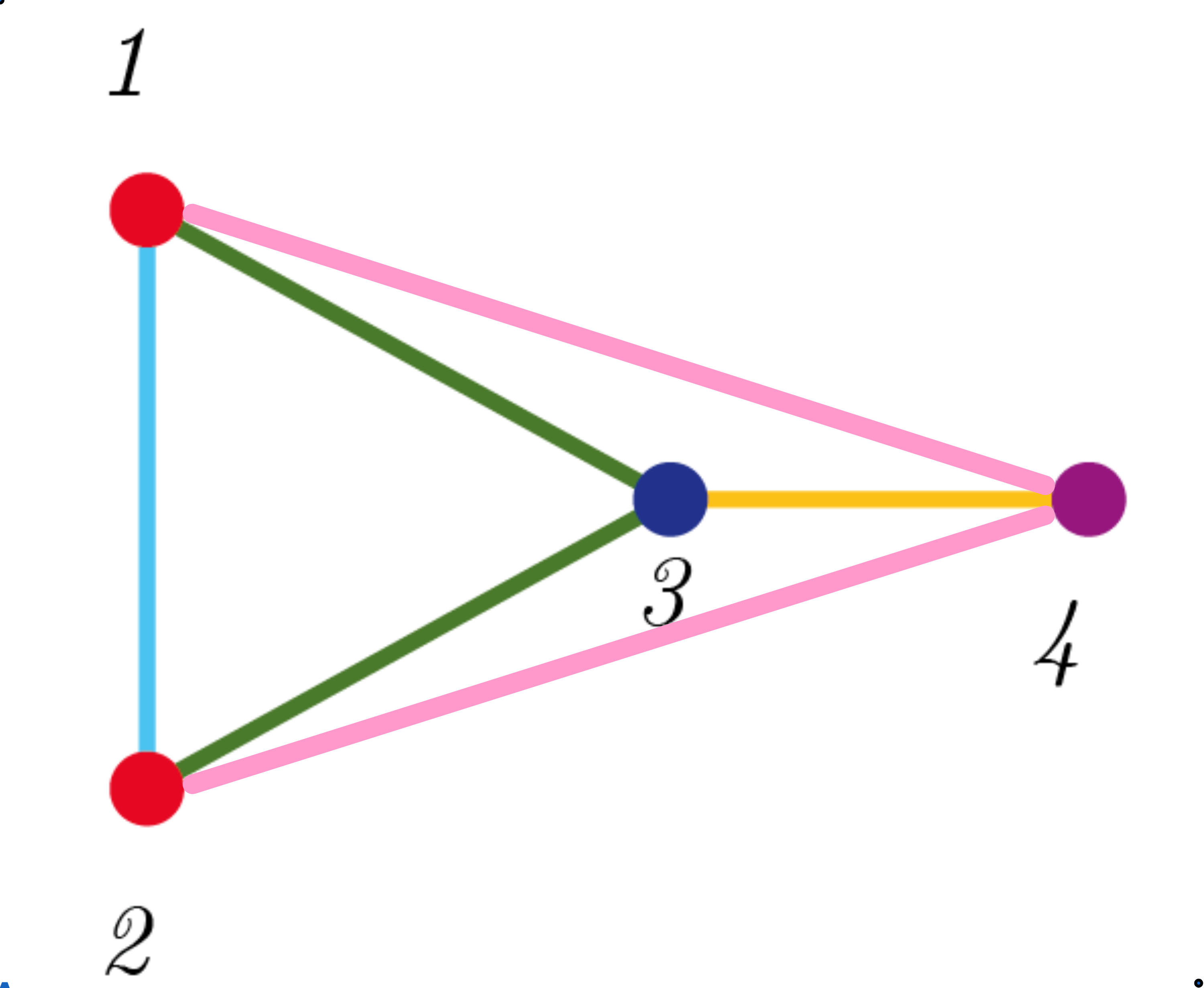}
    \caption{The RCOP block graph~$\cG$ from \Cref{fig:main} (on the left) and its completion~$\Gbar$ (on the right). These have vanishing ideals~\(I_{{\Gcal}} = \langle \sigma_{14} - \sigma_{24}, \sigma_{13} - \sigma_{23}, \sigma_{11} - \sigma_{22}, \sigma_{24}\sigma_{33}-\sigma_{23}\sigma_{34}\rangle\)  and~\(
I_{\overline{\Gcal}} = \langle \sigma_{14} - \sigma_{24}, \sigma_{13} - \sigma_{23}, \sigma_{11} - \sigma_{22} \rangle
\), respectively.}
    \label{fig:completion}
\end{figure}

\Cref{fig:completion} illustrates such a completion. Observe that the edges $\{1,4\}$ and $\{2,4\}$ of $\overline{\Gcal}$ have the same color since $\Lambda(1 \lra 4) = \Lambda(2 \lra 4)$.

\begin{theorem}\label{prop:completionRCOP}
The completion of an RCOP block graph is itself RCOP.
\end{theorem}

\begin{proof}
Let~$\overline{\cG}$ be the completion of the RCOP block graph~$\cG$, and let~$\Gamma(\Gcal)$ be the automorphism group associated to~$\cG$. We need to show that vertex and edge orbits of~$\Gamma(\Gcal)$ in~$\overline{\cG}$ are precisely its vertex and edge color classes, respectively. By definition of the completion, it is enough to prove that for any two edges of the same color in~$\overline{\cG}$, there is an automorphism~$\gamma$ in~$\Gamma(\Gcal)$ mapping one edge to the other. This is equivalent to proving that for any two combinatorially equivalent shortest paths in~$\cG$ there is an automorphism~$\gamma\in \Gamma(\Gcal)$ mapping one path to the other. This follows directly from
 \Cref{prop:ComboEquivalent} and \Cref{cor:completionRCOP}.
\end{proof}

Consider the \emph{Jordan algebra} structure on the set of $n \times n$ symmetric matrices $\S^n$ defined by
$X \bullet Y = \frac{1}{2}(XY + YX)$.
Note that this is the definition of a Jordan algebra given in \cite{bik2020jordan} with~$U$ equal to the identity matrix. For simplicity, we will refer to it as the Jordan algebra on~$\S^n$. Linear spaces of symmetric matrices that are  Jordan algebras are equal to their own inverse linear spaces. 

\begin{theorem}\label{thm:JordanAlgebras}\cite[Lemma 1]{jensen88}
Let~$\Lcal$ be a linear space of symmetric matrices. Let~$\Lcal^{-1}$ denote the Zariski closure of the set of all inverses of matrices in~$\Lcal$. Then~$\Lcal$ is a subalgebra of the Jordan algebra on~$\S^n$ if and only if~$\Lcal = \Lcal^{-1}$.
\end{theorem}

To check that a given linear space~$\Lcal$ forms a subalgebra of the Jordan algebra, one needs to prove that~$\Lcal$ is closed under the operation~$X \bullet Y$.
In fact, it suffices to show that for all~$X$ in~$\Lcal$, one has that~$X^2$ is in~$ \Lcal$ as well. 
This is a standard result in the study of Jordan subalgebras of symmetric matrices which we restate below for completeness.

\begin{prop}\label{prop:SquaresJordanAlg}
    Let~$\Lcal$ be a linear space of symmetric matrices. Suppose that for all~$X$ in~$\Lcal$,~$X^2$ also belongs to~$\Lcal$. Then~$\Lcal$ is a subalgebra of the Jordan algebra on~$\S^n$. 
\end{prop}

We can now prove that all RCOP complete graphs are generated by binomial linear relations determined by the coloring of the graph. 

\begin{theorem}\label{thm:rcop_complete}
Let~$\Gcal$ be a complete RCOP graph. Then~$\Lcal^{-1}_{\Gcal}$ is a Jordan algebra.
\end{theorem}

\begin{proof}
It is enough to prove that for any concentration matrix~$K$ in~$\Lcal_{\Gcal}$, the matrix~$K^2$ is also in~$\Lcal_{\Gcal}$. 
Denote the $(i,j)$ entry of $K$ by $K_{ij}$.
Note that since~${\Gcal}$ is complete, all entries of~$K$ are nonzero. Let~$r_i$ denote the~$i$-th row of~$K$;
since~$K$ is symmetric, this is also its~$i$-th column.
Thus the~$(i,j)$ entry of~$K^2$ is~$r_i \cdot r_j$.
In order to show that~$K^2$ is in~$\Lcal_{{\Gcal}}$, we must show that for each
$i, j\in [n]$ with~$\lambda(i) = \lambda(j)$, one has~$K^2_{ii} = K^2_{jj}$, and that for each pair of edges~$\{i,j\}$ and~$\{k,l\}$ in~$\cG$ with~$\lambda(\{i,j\}) = \lambda(\{k,l\})$,
one has~$K^2_{ij} = K^2_{kl}$.

Take vertices~$i$ and~$j$ with~$\lambda(i) = \lambda(j)$.
Since~$\Gcal$ is an RCOP graph, there exists an automorphism~$\gamma \in \Gamma(\Gcal)$
such that~$\gamma(i) = j$.
Thus we have 
\[
r_i \cdot r_i = \sum_{s = 1}^n m_{\lambda(\{i,s\})}^2 
= \sum_{s=1}^n m_{\lambda(\{\gamma(i),\gamma(s)\})}^2 
= \sum_{s=1}^n m_{\lambda(\{j,\gamma(s)\})}^2 
= \sum_{s=1}^n m_{\lambda(\{j,s\})}^2 
= r_j \cdot r_j.
\]
The second equality holds by definition of an RCOP graph.
The fourth equality holds because~$\gamma$ is a permutation on~$[n]$.

Now let~$\{i,j\}$ and~$\{k,l\}$ have~$\lambda(\{i,j\}) = \lambda(\{k,l\})$.
Then there is a~$\gamma$ in~$\Gamma(\Gcal)$ that maps one edge to the other, say,
that~$\gamma(i) = k$ and~$\gamma(j) = l$. 
Thus we have
\[
r_i \cdot r_j = \sum_{s=1}^n m_{\lambda(\{i,s\})} m_{\lambda(\{s,j\})} 
= \sum_{s=1}^n m_{\lambda(\{k, \gamma(s)\})} m_{\lambda(\{\gamma(s),l\})}  
= \sum_{s=1}^n m_{\lambda(\{k,s\})} m_{\lambda(\{s,l\})} 
= r_k \cdot r_l,
\]
as needed. The second equality holds since~$\gamma$ is an automorphism. The third equality holds because~$\gamma$ permutes the vertices of~$\Gcal$. The remainder of the theorem follows from \Cref{thm:JordanAlgebras} and \Cref{prop:SquaresJordanAlg}.
\end{proof}

Let us apply this result to completions of RCOP block graphs. By~ \Cref{prop:completionRCOP}, these are RCOP colorings of complete graphs. \Cref{thm:JordanAlgebras} and \Cref{thm:rcop_complete} together imply that their space of covariance matrices is equal to their linear space of concentration matrices. By definition of the completion, the latter is defined by the linear constraints~$K_{ij} =K_{i'j'}$ for each~$i,i',j,j'$ such that the shortest paths~$i\leftrightarrow j$ and~$i' \leftrightarrow j'$ are combinatorially equivalent. 
\begin{corollary}\label{cor:CompletionGenerators}
Let~$\overline{\Gcal}$ be the completion of the RCOP block graph~$\cG$. Then 
\begin{align}\label{eq:completeRCOP}
I_{\overline{\Gcal}} = \langle \sigma_{ij} - \sigma_{kl} \mid \Lambda(i\leftrightarrow j) = \Lambda(k \leftrightarrow l) \text{ in } \cG\rangle.
\end{align} 
\end{corollary}

We use~$\Bcal_{\Gbar}\coloneqq \{(i,j)-(k,l)\mid \Lambda(i\leftrightarrow j) = \Lambda(k \leftrightarrow l) \}$ to denote the Markov basis for~$I_{\Gbar}$ described in \Cref{eq:completeRCOP}.

\begin{example}\label{ex:completionIdeal} Consider the RCOP block graph~$\cG$ and its completion~$\Gbar$ in \Cref{fig:completion}. Then  \(
I_{\overline{\Gcal}} = \langle \sigma_{14} - \sigma_{24}, \sigma_{13} - \sigma_{23}, \sigma_{11} - \sigma_{22} \rangle.
\) The generators in this ideal are exactly the linear forms in~$I_{\cG}$ described in  \Cref{ex:main}.
\end{example}

\section{Another monomial map}\label{sec:AnotherMonomialMap}
We now introduce a new monomial map associated to an RCOP block graph~$\Gcal$. Its kernel coincides with~$\ker(\varphi_{\cG})$ from \Cref{def:shortestPathMapColored}.  This  second parametrization takes into account the color of \emph{all} vertices that lie on a  shortest path, and not only its endpoints.
This  small change will be of use in counting arguments when proving Theorem \ref{th:mainTheoremFormulation}. Let~$V(i\lra j)$ denote the set of vertices appearing in the shortest path~$i\lra j$.

\begin{definition}\label{def:newshortestPathMapColored}
The map~$\psi_\cG$ associated to a RCOP block graph~$\cG$ is the monomial map 
\begin{eqnarray*}
    \psi_{\cG}  :  \mathbb{R}[\sigma_{ij} \mid i,j \in [n], i \leq j] &\rightarrow& \mathbb{R}[t_{\ell}  \mid \ell\in\lambda(\cG)], \\
    \sigma_{ij} & \mapsto & \prod_{v \in V(i\leftrightarrow j)} t_{\lambda(v)} \prod_{e \in E(i\leftrightarrow j)} t_{\lambda(e)}.
\end{eqnarray*}
We refer to the exponent matrix for~$\varphi_{\cG}$  as~$B_\cG$.
\end{definition}
 The kernel of~$\psi_{\cG}$ is the toric ideal generated by all binomials of form
\begin{align}\label{eq:sp_ideal}  \sigma_{i_1j_1}\sigma_{i_2j_2}\dots \sigma_{i_dj_d}-\sigma_{k_1l_1}\sigma_{k_2l_2}\dots \sigma_{k_dl_d}, \end{align}
where the multisets~$L:=\{\{i_1\lra j_1,\dots, i_d\lra j_d\}\}$ and~$R:=\{\{k_1\lra l_1,\dots, k_d\lra l_d\}\}$ have the same multiset of vertex and edge colors. 

\begin{theorem}\label{cor:counting}
Let~$\cG$ be an RCOP block graph. Then~$\ker(\psi_\cG)=\ker(\varphi_\cG)$. In particular, if~$f=f_1-f_2$ is a binomial in~$\ker (\varphi_\cG)$, then the multisets of vertex colors in the shortest paths used by~$f_1$ and~$f_2$, respectively, are the same. 
\end{theorem}
\begin{proof}
It is enough to prove that the matrices~$A_{\cG}$ and~$B_\cG$ have the same rowspan. We will do this by showing that~$B_\cG$ is obtained by applying elementary row operations to~$A_\cG$. This completes the proof as elementary row operations are invertible.

Let~$a_{\ell}$  be the row of~$A_\cG$ corresponding to the color~$\ell$ in~$\lambda(\cG)$.  Similarly we have~$b_{\ell}$  for the rows of~$B_\cG$. 
It is clear by the definitions of the two maps that~$a_{\epsilon}=b_{\epsilon}$ for any edge color~$\epsilon$ in~$ \lambda(\cG)$. We will show that for any vertex color~$\nu$ in~$\lambda(\cG)$,
\begin{align} \label{eq:rowspan}
b_\nu&=& \dfrac{1}{2} \left( a_\nu + \sum_{\substack{ \text{edge }e: \text{ one end of }e \\ \text{ has color }\nu} }a_{\lambda(e)} +2\sum_{\substack{\text{edge }e: \text{ both ends of }e \\ \text{ have color }\nu } }a_{\lambda(e)}\right). 
\end{align}
Let us take two vertices~$i$ and~$j$ and consider all possible 
scenarios on the
position of vertices of color~$\nu$ in the shortest path~$i\lra j$.  Notice that due to \Cref{lm:twovertices}, there cannot be more than two vertices of color~$\nu$ in~$i\lra j$.  

\emph{Case 1}: Suppose that the color~$\nu$ appears only once in~$i \lra j$. This means that the $\{i,j\}$ entry of~$b_\nu$ is $b_\nu(\{i,j\})=1$. 

First we consider the case where either~$i$ or~$j$ has color~$\nu$. If~$i=j$, the path~$i\lra j$ is just a vertex. So,~$a_{\nu}(\{i,j\})=2$ and  the left-hand side of~\Cref{eq:rowspan}  is~$(2+0+0)/2=1$. 
If~$i\neq j$, without loss of generality we may assume that~$\lambda(i) = \nu$. In this case~$a_{\nu}(\{i,j\})=1$.  The  only edge in~$i\lra j$ having at least one endpoint of color~$\nu$ is the one containing vertex~$i$. Hence, the sum in the left-hand side of~(\ref{eq:rowspan}) is exactly~$(1+1+0)/2=1$. 

Now suppose that the color~$\nu$ appears at one of the internal nodes~$i'$ of~$i\lra j$. In this case,~$a_\nu(\{i,j\})=0$. The only edges in~$i\lra j$ having at least one endpoint of color~$\nu$ are the two edges that contain~$i'$. Note that the other endpoints of these two edges do not have color~$\nu$, since otherwise this color appears more than once in~$i\lra j$. Hence, the right-hand side of~(\ref{eq:rowspan}) is~$ (0+2+0)/2=1$.

\emph{Case 2}: Suppose that the shortest path~$i\lra j$ contains two vertices~$i'$ and~$j'$ of color~$\nu$. Here we have~$b_{\nu}(\{i,j\})=2$.  

If~$i'$ and~$j'$  are adjacent to each other,  they are the endpoints of edge~$e = \{i',j'\}$ in~$i\lra j$. If~$i = i'$ and~$j = j'$, then both~$a_{\nu}(\{i,j\})$  and~$b_{\nu}(\{i,j\})$ are equal to~$2$. 
If~$i\lra j$ contains more than one edge and~$i = i'$  (similarly for~$j=j'$), then~$a_{\nu}(\{i,j\})$ is~$1$. The edge following~$e$  in~$i\lra j$ has exactly one endpoint of color~$\nu$. Hence, the right-hand side of~(\ref{eq:rowspan}) is~$(1+1+2)/2=2$. Finally, if~$e$ is an internal edge of~$i\lra j$, 
then the two other edges adjacent to~$e$ in~$i\lra j$ contain exactly one vertex 
of color~$\nu$. In this case,~$a_{\nu}(\{i,j\})=0$ and the right-hand side of ~(\ref{eq:rowspan}) is~$(0+2 + 2)/2=2$.

Suppose now that~$i'$ and~$j'$ are not adjacent to each other. If~$i = i'$ and~$j = j'$, then~$a_{\nu}(\{i,j\})=2$ and the right-hand side of~(\ref{eq:rowspan}) is~$(2+2 + 0)/2=2$. If~$i = i'$ but~$j \neq j'$ (similarly for~$i\neq i'$ and~$j= j'$), then there are in total three edges having exactly one endpoint of color~$\nu$; the one containing~$i$ and the two containing~$j'$. Hence, we have~$a_{\nu}(\{i,j\})=1$ and the right-hand side of~(\ref{eq:rowspan}) is~$(1+3 +0)/2=2$. Lastly, if~$i'$ and~$j'$ are both internal vertices of~$i\lra j$, then there are in total four edges having exactly one endpoint of color~$\nu$. So,~$a_{\nu}(\{i,j\})=0$ and the right-hand side of~(\ref{eq:rowspan}) is~$(0+4 +0)/2=2$. 

We have exhausted all possible scenarios. Thus~$B_\cG$ is obtained by applying elementary row operations to~$A_\cG$, and hence the ideals~$\ker(\varphi_\cG)$ and~$\ker(\psi_\cG)$ are equal.
\end{proof}

\section{Markov bases of RCOP block graphs}\label{sec:Main}

Now we are ready to prove the main result of this paper, namely, that $I_\cG$ is equal to the toric ideal $I_G+I_{\overline{\Gcal}}$. The main idea is to show that $\ker(\varphi_{\cG})=I_\cG+I_{\overline{\Gcal}}\subseteq I_{\cG}$.  Then we complete the proof by showing that the prime ideals $\ker(\varphi_{\cG})$ and $I_\cG$ have the same dimension.

\begin{theorem}\label{th:mainTheoremFormulation}
Let~$\cG$ be an RCOP block graph with underlying graph~$G$ and completion~$\Gbar$. Then the union of Markov bases~$\mathcal{B}_{G}$ in (\ref{eq:gensG}) and~$\mathcal{B}_{\overline{\cG}}$ in (\ref{eq:completeRCOP}) is a Markov basis for~$\ker(\varphi_{\cG})$. 
\end{theorem}

In order to prove this theorem, we choose an arbitrary  binomial in $\ker(\phi_{\Gcal})$ and show that we can apply a sequence of Markov moves from $\Bcal_G \cup \Bcal_{\overline{\Gcal}}$ to reduce this binomial to $0$. We establish the following notation which we will use for the remainder of the section.

\begin{hyp}\label{hyp1}
Let $f = \sigma_{i_1j_1}\cdots \sigma_{i_dj_d} - \sigma_{k_1l_1}\cdots \sigma_{k_dl_d}$ be a binomial in the homogeneous toric ideal $\ker(\varphi_{\Gcal})$. Let~$L$ denote the multiset of the shortest paths~$i_1\leftrightarrow j_1$,~$\ldots$,~$i_d\leftrightarrow j_d$ appearing in the left-hand side monomial in $f$. Similarly, let ~$R$ be the multiset of paths~$k_1\leftrightarrow l_1$,~$\ldots$,~$k_d\leftrightarrow l_d$ on the right hand-side of $f$. We further assume that $i_1 = k_1$ but $j_1 \neq l_1$. To avoid heavy notation, assume that the vertices of $k_1 \lra l_1$ are $k_1 = 1, 2, \dots, l_1$ appearing in this order. Let $c$ be the largest vertex in $k_1 \lra l_1$ that also lies on $i_1 \lra j_1$ so that $k_1 \lra c = i_1 \lra c$; this is consistent with the labeling in \Cref{fig:OriginalPaths}. Assume that $l_1 > c$ so that $c+1$ is the vertex on $k_1 \lra l_1$ following $c$. Let $a$ be the vertex on $i_1 \lra j_1$ following $c$, if it exists. Further assume that $\lambda(\{c,c+1\}) \neq \lambda(\{c,a\})$ if both edges exist.
\end{hyp}

We aim to show that there exists a sequence of Markov moves that we can apply to $f$ to obtain a path $i_1 \lra j_1'$ that matches with $k_1 \lra l_1$ on one further edge. The following hypothesis asserts the existence of a second path in $L$ that we will use to construct these Markov moves; we prove that this hypothesis holds in the proof of \Cref{th:mainTheoremFormulation}.

\begin{hyp}\label{hyp2}
Let $i_2 \lra j_2 \in L$ contain edge $\{c, c+1\}$ appearing in that order. If $i_2 \neq c$, let $b$ denote the vertex in $i_2 \lra j_2$ that appears immediately before $c$.
\end{hyp}

We can now examine the Markov moves that will replace $i_1 \lra j_1$ with a path $i_1 \lra j_1'$ that matches $k_1 \lra l_1$ up to one more edge. There are two main cases depending upon whether $a$ and $b$ are adjacent or not. We first consider the case where $a$ and $b$ are not adjacent.

\begin{lemma}\label{lm:nonadjacent}
Adopt \Cref{hyp1} and \Cref{hyp2}. Suppose that either $a$ or $b$ do not exist, or that they both exist but $a$ and $b$ are not adjacent. Then there exists a Markov move in $\Bcal_G$ 
that we may apply to $(i_1,j_1)\dots(i_d,j_d)$ in order to replace $(i_1,j_1)$ with $(i_1,j_1')$ such that $i_1 \lra j_1'$ passes through $c+1$.
\end{lemma}

\begin{figure}[ht]
    \centering
    \includegraphics[scale=.6]{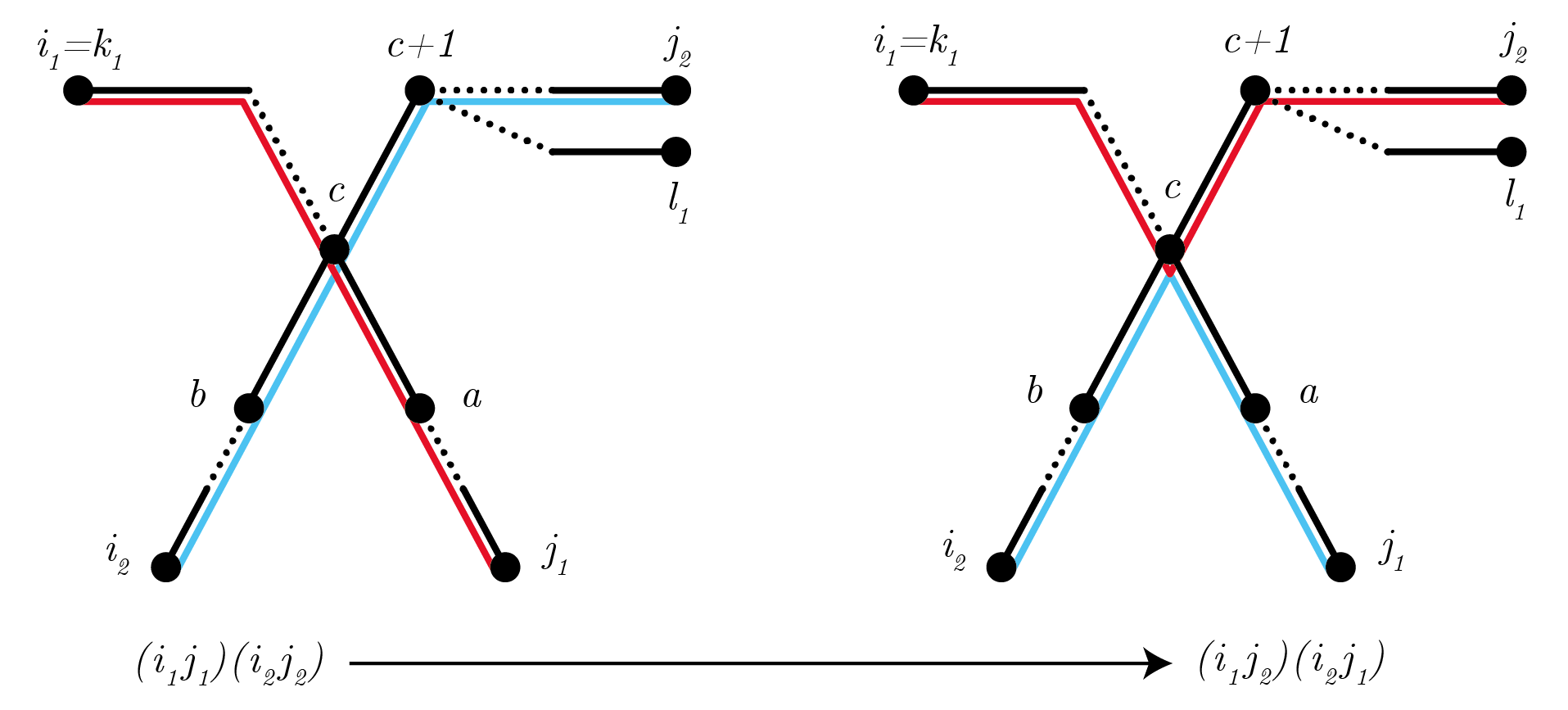}
    \caption{Markov move when~$[a,b]$ is not an edge in~$G$.}
    \label{fig:NoEdgeab}
\end{figure}

\begin{proof}
We claim that the basis~$\mathcal{B}_{G}$ contains the quadratic Markov move \(b = (i_1,j_1)(i_2,j_2)-(i_1,j_2)(i_2,j_1)\) which yields the desired result. This move is visualized in \Cref{fig:NoEdgeab}.

Indeed, as $i_1 = k_1$, we have that $i_1 \lra c+1$ contains the edge $\{c,c+1\}$. Similarly since $b$ and $a$ are not adjacent, the shortest path $i_2 \lra a = i_1 \lra \{c,a\}$. Thus by \Cref{lm:union}, the multiset of edges in $i_1 \lra j_1$ and $i_2 \lra j_2$ is the same as that of $i_1 \lra j_2$ and $i_2 \lra j_1$. Thus $b \in \Bcal_G$ and $(i_1,j_2)$ contains the edge ${c,c+1}$, as needed.
\end{proof}

We now have two lemmas regarding the case where $a$ and $b$ are adjacent. They each deal with two separate cases regarding the structure of the cliques that meet at $c$.

\begin{lemma}\label{lem:severalcliques}
Adopt \Cref{hyp1} and \Cref{hyp2}. Suppose that $a$ and $b$ exist and that $a,b$ and $c$ belong to a maximal clique $\cC$. Suppose further that there is another clique $\cC'$ that is isomorphic to $\cC$ and contains $c$. Then there exists a sequence of Markov moves in $\Bcal_G \cup \Bcal_{\overline{\Gcal}}$ 
that we may apply to $(i_1,j_1)\dots(i_d,j_d)$  to replace $(i_1,j_1)$ with $(i_1,j_1')$ such that $i_1 \lra j_1'$ passes through $c+1$.
\end{lemma}

\begin{figure}[ht]
    \centering
    \includegraphics[scale=.7]{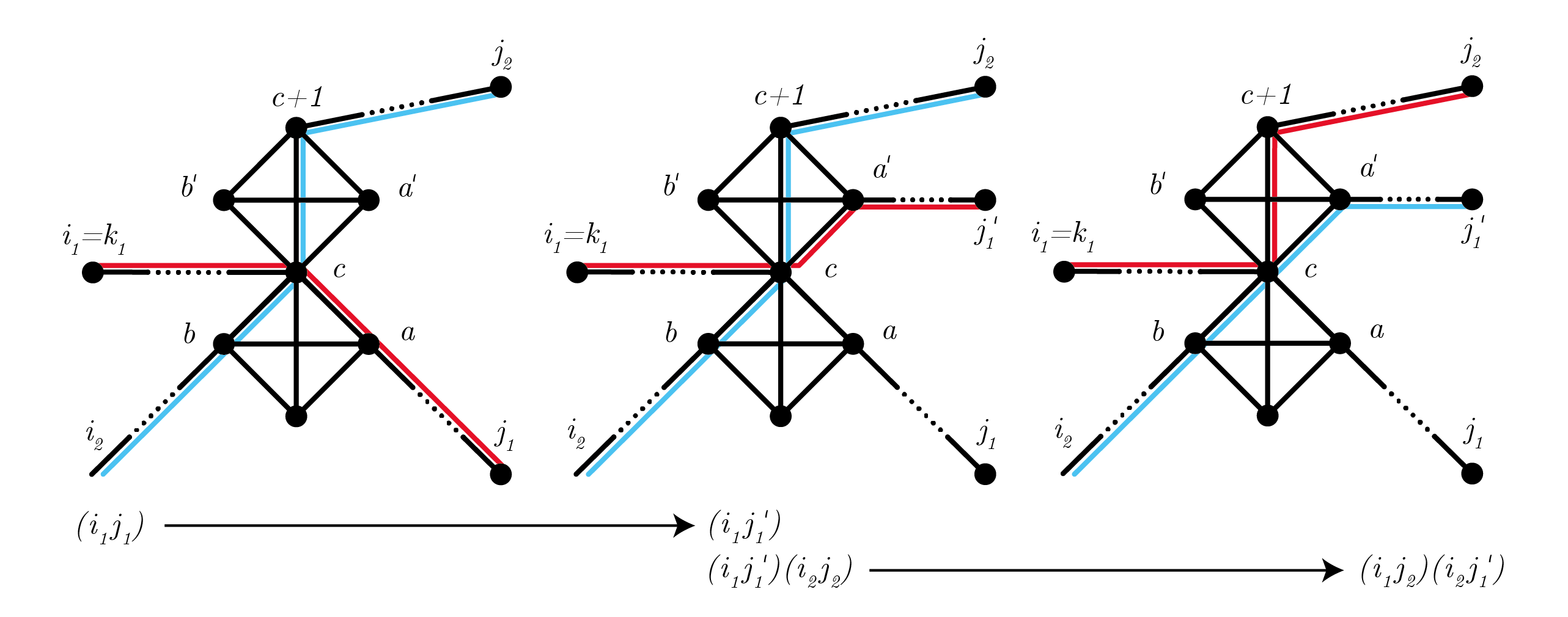}
    \caption{Markov move for Case 1 in the proof of \Cref{lem:severalcliques}.}
    \label{fig:Case1a}
\end{figure}

\begin{figure}[ht]
    \centering
    \includegraphics[scale = .7]{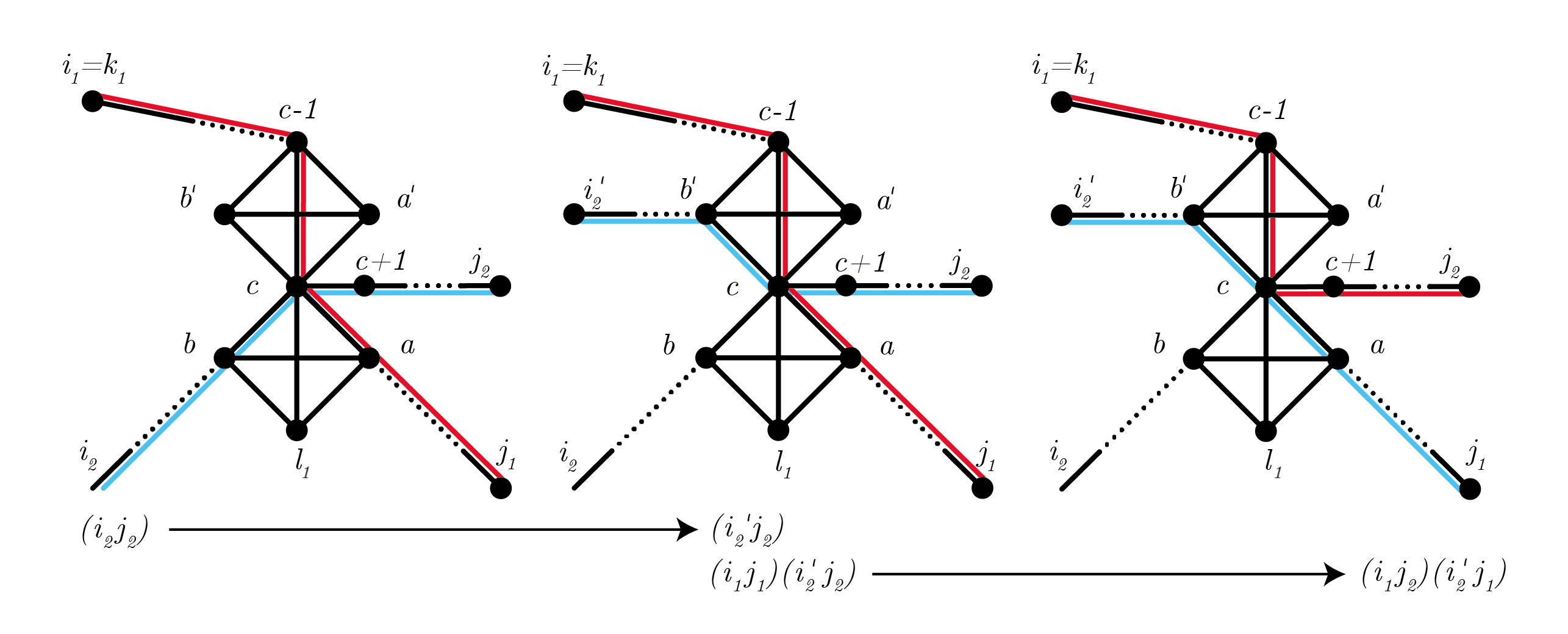}
    \caption{Markov moves for Case 2. in the proof of \Cref{lem:severalcliques}.}
    \label{fig:Case1b}
\end{figure}

\begin{proof}
The vertices~$c-1$ and~$c+1$ cannot both be contained in $\cC'$ since they are not adjacent. Let~$a', b'$ be two vertices in~$ \cC'$ such that~$\lambda(\{c,a\}) = \lambda(\{c,a'\})$ and~$\lambda(\{c,b\}) = \lambda(\{c,b'\})$. This implies that~$\lambda(a) = \lambda(a')$ and~$\lambda(b) = \lambda(b')$. 
We have two cases based on the position of~$c+1$.

\emph{Case 1.} Suppose that~$c+1 \in \cC'$. Let~$d \in \cC$ such that~$\lambda(\{c,d\}) = \lambda(\{c,c+1\})$. Thus~$\lambda(d) = \lambda(c+1)$. 
By \Cref{lm:alphabeta}\ref{lm:beta}, there exists a linear move~$(i_1, j_1) - (i_1, j_1')$ that fixes~$i_1 \lra c$ and sends~$\{c,a\}$ to~$\{c,a'\}$. We apply this move.
Since~$\cC'$ is a clique, we have $i_1 \lra j_2 = i_1 \lra c\lra c+1 \lra j_2$ and~$i_2 \lra j_1' = i_1 \lra c\lra a' \lra j_1'$. 
Thus~$(i_1, j_1') (i_2, j_2) - (i_1, j_2) (i_2, j_1')$ is a move in~$\mathcal{B}_G$. We apply this move so that~$i_1\lra  j_2$ matches~$k_1\lra l_1$ up to one more edge, as needed. This is pictured in \Cref{fig:Case1a}.

\emph{Case 2:} Now suppose that~$c+1 \not\in \cC'$.
Then~$c+1$ lies in a~$c$-component different from~$\cC$ and~$\cC'$. Thus, by \Cref{lm:alphabeta}\ref{lm:beta} we know that there exists a linear move~$(i_2, j_2) - (i_2', j_2)$ that fixes~$c \lra j_2$ and sends vertex~$b$ to~$b'$. 
Since~$\cC'$ is a clique,~$(i_1, j_1) (i_2', j_2) - (i_1, j_2) (i_2', j_1)$ is a Markov move in~$\mathcal{B}_G$. We apply this move so that~$i_1 \lra j_2$ matches~$k_1\lra l_1$ in one more edge, as needed.  Note that these moves can be applied whether~$c-1 \in \cC'$ (as pictured in \Cref{fig:Case1b}), or not. 
\end{proof}

Finally, we consider the most technical case, in which there is no other clique isomorphic to $\cC$ that is adjacent to $c$.

\begin{lemma}\label{lem:oneclique}
Adopt \Cref{hyp1} and \Cref{hyp2}. Suppose that $a$ and $b$ exist and that $a,b$ and $c$ belong to a maximal clique $\cC$. Suppose further that $\cC$ is the only clique in its isomorphism class that contains $c$. Then there exists a sequence of Markov moves in $\Bcal_G \cup \Bcal_{\overline{\Gcal}}$ 
that we may apply to $(i_1,j_1)\dots(i_d,j_d)$ to replace $(i_1,j_1)$ with $(i_1,j_1')$ such that $i_1 \lra j_1'$ passes through $c+1$.
\end{lemma}

\begin{proof}
Denote by $\lambda(\cC)$ the set of colors of edges in $\cC$.
We will show that either (i) there exists a path in~$L$ that contains a vertex with color~$\lambda(c)$ whose adjacent edges are not colored with any colors from~$\lambda(\cC)$, or (ii) there is an automorphism in~$\Gamma(\cG)$ that allows us to apply one of the previous lemmas. 

We first make an important observation which will be used to prove this case.
Let~$c'$ be a vertex with~$\lambda(c') = \lambda(c)$. Let~$i \lra j$ be a shortest path containing~$c'$ and let~$d$ and $e$ be the vertices in~$i \lra j$ appearing directly before and after~$c'$.
We claim that at most one of colors~$\lambda(\{d,c'\})$ and~$\lambda(\{c',e\})$ can belong to~$\Lambda(\cC)$. 
Suppose for contradiction that both~$\lambda(\{d,c'\})$ and~$\lambda(\{c',e\})$ are present in~$\lambda(\cC)$. Since~$\{d,c'\}$ and~$\{c',e\}$ lie in a shortest path we have that~$d$ and~$e$ cannot be in a single clique. Let~$\cC_1$ and~$\cC_2$ be the two distinct cliques containing~$\{d,c'\}$ and~$\{c',e\}$ respectively. As~$\cC$ is the only clique in the isomorphism class containing~$c$, 
by \Cref{lm:iso_cliques} we know that both~$\cC_1$ and~$\cC_2$ have to be isomorphic to~$\cC$. Under these conditions, the two vertices~$c$ and~$c'$ of the same color have non isomorphic neighborhoods which contradicts  \Cref{prop:RCOPneighb}. 

We introduce the following useful notation. Let~$\lambda^{cc}(\cC)$ denote the set of all edge colors in~$\cC$ that connect two vertices of color~$\lambda(c)$. Let~$\lambda^c(\cC)$ denote the set of all edge colors in~$\cC$ that connect a vertex of color~$\lambda(c)$ to a vertex whose color is \emph{not}~$\lambda(c)$. Let~$E(L)$ and~$E(R)$ be the multiset of all edges in paths in~$L$ and~$R$, respectively. Similarly, let~$V(L)$ and~$V(R)$ denote the multiset of all vertices in paths in~$L$ and~$R$. Finally, let~$U_{c,\cC}(L)$ denote the multiset of all vertices in paths in~$L$ with color~$\lambda(c)$ whose adjacent edges in their respective paths do \emph{not} have colors in~$\lambda(\cC)$. Define~$U_{c,\cC}(R)$ similarly. 

Let~$c'$ be a vertex of color~$\lambda(c)$ and let~$i \lra j$ be a path in~$L$ that contains it. By \Cref{prop:RCOPneighb} vertex~$c'$ is contained in exactly one clique isomorphic to~$\cC$. Hence, exactly one of the following three conditions must hold: either~$c'$ is adjacent in~$i\lra j$ to one edge whose color belongs to~$\lambda^{cc}(\cC)$, or~$c'$ is adjacent in~$i \lra j$ to one edge whose color belongs to~$\lambda^c(\cC)$, or~$c'$ belongs to~$U_{c,\cC}(L)$. Note that the vertex~$c'$ cannot be simultaneously adjacent to edges with colors in~$\lambda^{cc}(\cC)$ and~$\lambda^{c}(\cC)$ in~$i\lra j$ as the clique~$\cC$ is the only clique in its isomorphism class that contains~$c$.  Hence the number of vertices in~$L$ of color~$\lambda(c)$ is equal to
\[
2\#\{\{e\in E(L)~|~\lambda(e)\in \Lambda^{cc}(\cC)\}\} + \#\{\{e\in E(L)~|~\lambda(e)\in \Lambda^{c}(\cC)\}\} + \#U_{c,\cC}(L).
\]

Similarly, the number of vertices in~$R$ of color~$\lambda(c)$ is equal to
\[
2\#\{\{e\in E(R)~|~\lambda(e)\in \lambda^{cc}(\cC)\}\} + \#\{\{e\in E(R)~|~\lambda(e)\in \lambda^{c}(\cC)\}\} + \#U_{c,\cC}(R).
\]

By definition of the shortest path map,
we have that~$\#(\lambda^{cc}(\cC) \cap E(L)) = \#(\lambda^{cc}(\cC) \cap E(R))$ and~$\#(\lambda^c(\cC) \cap E(L)) = \#(\lambda^c(\cC) \cap E(R))$. Moreover, by \Cref{cor:counting}, the numbers of vertices of color~$\lambda(c)$ in~$L$ and~$R$ are equal. Hence, it must be the case that~$\#U_{c,\cC}(L) = \#U_{c,\cC}(R)$. Moreover, we know  by construction that none of the edges in~$\cC$ can lie in~$k_1\lra l_1$. This implies that in~$k_1 \lra l_1$,~$c$ is not adjacent to an edge in~$\cC$. Hence,~$\#U_{c,\cC}(L) = \#U_{c,\cC}(R) \geq 1$.
So there exists a shortest path~$i \lra j$ in~$L$ that contains a vertex of color~$\lambda(c)$ which is not adjacent to an edge in~$i \lra j$ whose color is in~$\lambda(\cC)$. Let~$c'$ denote this vertex.

There are three scenarios for the path $i\lra j$: either~$i \lra j = i_1 \lra j_1$,~$i \lra j = i_2 \lra j_2$ or~$i \lra j = i_m \lra j_m$ for~$m \geq 3$. 

\emph{Case 1:} Suppose that~$i \lra j = i_1 \lra j_1$. Without loss of generality, we may assume that~$c'$ does not appear in~$i_1 \lra c$. This is because $i_1 = k_1$, so $c'$ is also a vertex in $k_1 \lra c$. Thus $\#U_{c,\cC}(L) = \#U_{c,\cC}(R) \geq 2$ implies that there is another vertex in~$L$ of color $\lambda(c)$ with the desired property. 
Hence~$c'$ belongs to~$c \lra j_1$. The subpath~$c \lra c'$ is symmetric by \Cref{lm: symmetry}. Hence~$c'$ is adjacent to an edge of color~$\lambda([c,a])$ in~$i_1 \lra j_1$, which is a contradiction. Thus we cannot, in fact, have that~$i \lra j$ is $i_1 \lra j_1$.

\emph{Case 2:} Suppose that~$i \lra j = i_2 \lra j_2$. The vertex~$c'$ cannot lie in~$i_2 \lra c$. Indeed, since the subpath~$c' \lra c$ is symmetric, this would imply that~$c'$ is adjacent  to an edge of color~$\lambda(\{b,c\})$ in~$i_2 \lra j_2$, which is a contradiction. So,~$c'$ must lie in~$c \lra j_2$. Let~$(c+1)'$ denote the vertex adjacent to~$c'$ that lies on the path~$c \lra c'$; by the symmetry of this path, we have~$\lambda(\{(c+1)',c'\}) = \lambda(\{c,c+1\})$. Let~$d$ denote the vertex after~$c'$ in~$i_2 \lra j_2$ if it exists. The cases where~$(c+1)' = c$ and~$(c+1)' \neq c$ are slightly different.

First suppose that~$(c+1)' \neq c$. Then, by  \Cref{lm:twovertices},~$\lambda((c+1)')$ must be different from~$\lambda(c)$. There exists an automorphism~$\gamma \in \Gamma(\cG)$ that maps~$\{c', (c+1)'\}$ to~$\{c,c+1\}$. This automorphism must send~$c'$ to~$c$ and~$(c+1)'$ to~$c+1$. We may replace~$i_2 \lra j_2$ with~$\gamma(j_2) \lra \gamma(i_2)$ in~$L$. If~$d$ exists, then~$\lambda(\{c',d\}) \not\in \Lambda(\cC)$, so~$\gamma(d) \not\in \cC$. Hence, we are now in the case of \Cref{lm:nonadjacent}. We may apply the Markov move~$(i_1,j_1)(\gamma(j_2),\gamma(i_2))- (i_1, \gamma(i_2))(\gamma(j_2), j_1)$ from $\mathcal{B}_G$, and~$i_1 \lra \gamma(i_2)$ matches~$k_1\lra l_1$ up to one more edge, as needed. 

Now suppose that~$(c+1)' = c$. In other words, we have that~$c' = c+1$. There exists an automorphism~$\gamma \in \Gamma(\cG)$ that sends~$c+1$ to~$c$. By \Cref{lm:alphabeta}\ref{lm:beta}, there exists an automorphism~$\alpha \in \Gamma(\cG)$ that fixes~$\gamma(j_2) \lra c$ and maps~$\gamma(c)$ to~$c+1$. Hence the path~$\alpha(\gamma(j_2)) \lra \alpha(\gamma(i_2))$ has~$\{c,c+1\}$ as an edge and~$c$ is not adjacent to a vertex of~$\cC$ in this path. Thus we may replace~$i_2\lra j_2$ with~$\alpha(\gamma(j_2))\lra \alpha(\gamma(i_2))$ in~$L$ and then apply the Markov move~$(i_1,j_2)(\alpha(\gamma(j_2)),\alpha(\gamma(i_2))) - (i_1, \alpha(\gamma(i_2)))(\alpha(\gamma(j_2)),j_1)$ from~$\mathcal{B}_G$. Thus,~$i_1\lra  \alpha(\gamma(i_2))$ matches~$k_1\lra l_1$ up to one more edge, as needed.

\emph{Case 3:} The rest of the proof is concerned with the case where~$i \lra j = i_m \lra j_m$ for~$m \geq 3$. First, without loss of generality, we assume that~$i \lra j = i_3 \lra j_3$. Moreover, we may assume that~$c' = c$. Indeed, there exists an automorphism~$\gamma \in \Gamma(\cG)$ that sends~$c'$ to~$c$, so we may simply replace~$i_3 \lra j_3$ with~$\gamma(i_3) \lra \gamma(j_3)$.

Let~$d$ and~$e$ denote those vertices that appear in~$i_3 \lra j_3$ immediately before and after~$c$, respectively. Note that~$\{a,d\}$,~$\{a,e\}$,~$\{b,d\}$ and~$\{b,e\}$ are not edges in~$\Gcal$, as this would imply that~$d$ or~$e$ belong to~$\cC$. Thus we have the following relationships between shortest paths:
\begin{align*}
    i_2 \lra j_3 & = i_2 \lra c \lra j_3,  \\
    i_3 \lra j_1 &= i_3 \lra c \lra j_1, \\
    i_2 \lra i_3 &= i_2 \lra c \lra i_3, \text{ and}\\
    j_1 \lra j_3 &= j_1 \lra c \lra j_3.
 \end{align*} 
We now have two subcases based on whether~$\{d,c+1\}$ is an edge in~$\Gcal$.
 
\emph{Case 3a:} Suppose that~$\{d,c+1\}$ is not an edge in~$\Gcal$. In particular, this implies that
\(
i_3 \lra j_2 = i_3 \lra c \lra j_2.
\)
Thus, we may apply to~$L$ the moves~$(i_2, j_2) (i_3, j_3) - (i_2, j_3) (i_3, j_2)$ and then~$(i_1, j_1)(i_3, j_2) - (i_1, j_2) (i_3, j_1)$ from~$\mathcal{B}_G$.
This transforms~$L$ into
 \[
L= \{i_1\lra j_2, i_2\lra j_3, i_3\lra j_1,  i_m\lra j_m, \text{for } m=4,5,\dots,d\}.
 \]
 Now~$i_1 \lra j_2$ matches~$k_1 \lra l_1$ up to one more edge, as needed. This is illustrated in \Cref{fig:Case2a}.
 
 \begin{figure}[ht]
     \centering
     \includegraphics[scale=.7]{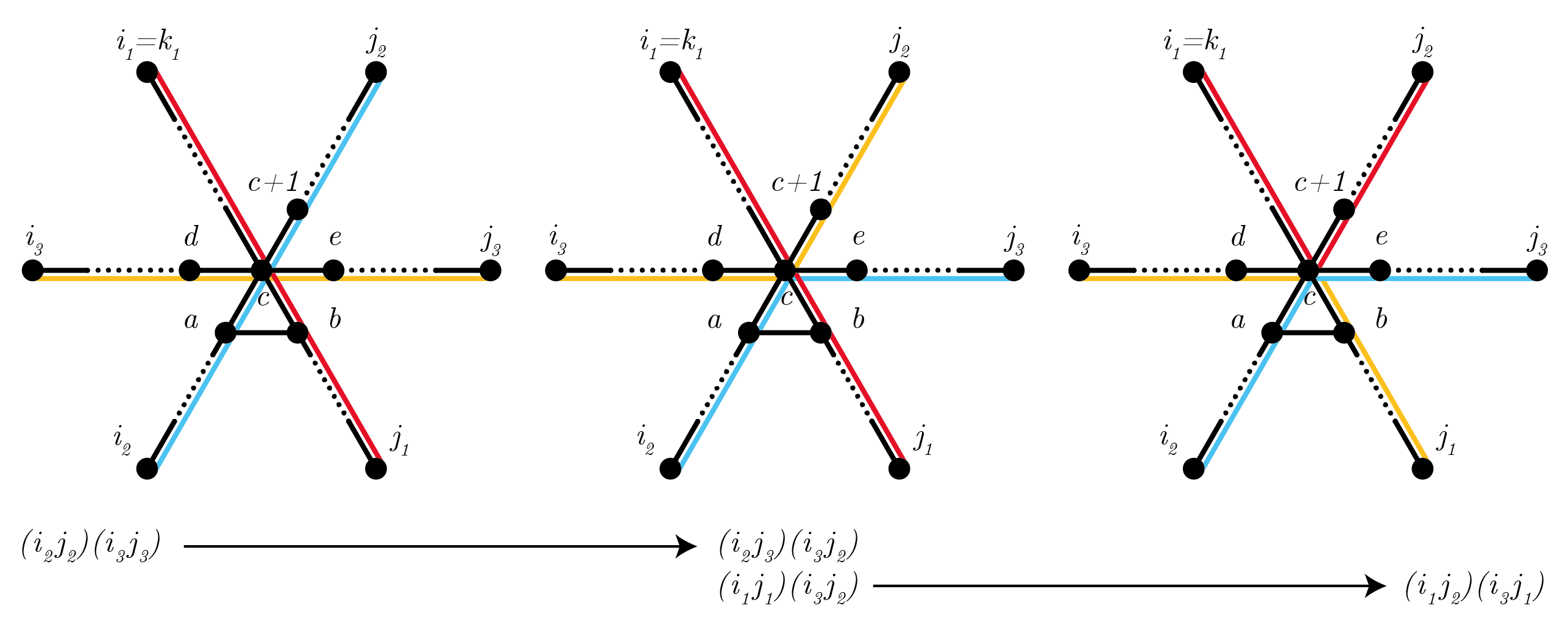}
     \caption{Markov moves for case 3a in the proof of \Cref{th:mainTheoremFormulation}.}
     \label{fig:Case2a}
 \end{figure}
 
\begin{figure}[ht]
     \centering
     \includegraphics[scale=.7]{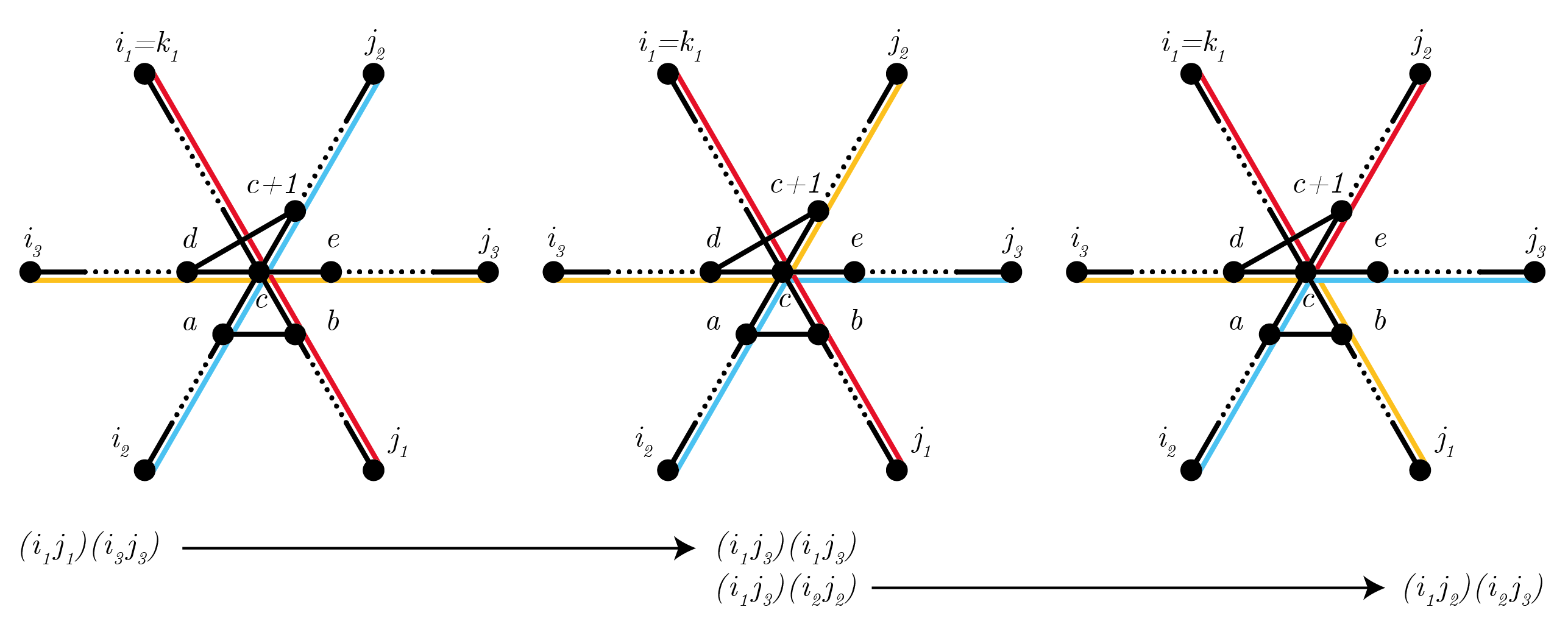}
     \caption{Markov moves for Case 3b in the proof of \Cref{th:mainTheoremFormulation}.}
     \label{fig:Case2b}
 \end{figure}

\emph{Case 3b:} Now suppose that~$\{d,c+1\}$ is an edge of~$\Gcal$. This implies that~$\{d,c-1\}$ is not an edge of~$\Gcal$ because if such an edge existed, we would need to have a clique that contains~$d, c-1, c$ and~$c+1$. This contradicts that~$c-1$ and~$c+1$ are not adjacent. Thus we have that 
\(
i_1 \lra i_3 = i_1 \lra c \lra i_3.
\)
So, we may apply in sequence the moves~$(i_1, j_1) (i_3, j_3) - (i_1, i_3) (j_1, j_3)$ and~$(i_1, i_3)(i_2, j_2) - (i_1, j_2) (i_2, i_3)$ form $\mathcal{B}_G$. 
This transforms~$L$ into
\[
L= \{i_1\lra j_2, i_2\lra i_3, j_1\lra j_3,  i_m\lra j_m, \text{for } m=4,5,\dots,d\}.
 \]
 Now~$i_1 \lra j_2$ matches~$k_1 \lra l_1$ up to one more edge, as needed. This is illustrated in \Cref{fig:Case2b}.
 \end{proof}
 
We have now amassed all of the results necessary to prove \Cref{th:mainTheoremFormulation}.

\begin{proof} 
Let $f = \sigma_{i_1j_1}\cdots \sigma_{i_dj_d} - \sigma_{k_1l_1}\cdots \sigma_{k_dl_d}$ be a binomial in $\ker(\varphi_{\Gcal}).$ Since the toric ideal $\ker(\varphi_{\Gcal})$ is homogeneous, all binomials in it have this form.
Let $L$ and $R$ denote the multisets of the shortest paths appearing in the left-hand side and right-hand side monomials in $f$, respectively, as in the previous lemmas.
We show by induction on the number $d$ of paths in~$L$ and~$R$ that $f$ can be reduced to the $0$ polynomial via Markov moves in $\Bcal_G \cup \Bcal_{\overline{\Gcal}}$. Notice that $d$ is also the degree of the binomial $f$.

 If~$d=1$, then the two paths~$i_1 \leftrightarrow j_1$ and~$k_1 \lra l_1$ must be  combinatorially equivalent. Thus by \Cref{cor:CompletionGenerators}, there exists an automorphism in~$\Gamma(\cG)$ mapping one path to the other. 
Take now~$d>1$. For path~$i_1\leftrightarrow j_1$ in~$L$, there is at least one path in~$R$ with at least one endpoint of  color~$\lambda(i_1)$. Without loss of generality, assume that this path is~$k_1\leftrightarrow l_1$ and~$\lambda(k_1) = \lambda(i_1)$.
To avoid heavy notation, assume further that vertices of the second path are~$k_1=1,2,\ldots,l_1$, appearing in this order.  
Take the largest vertex~$c$ in~$k_1\leftrightarrow l_1$ for which the coloring in the subpath~$k_1\leftrightarrow c$ matches component-wise with the coloring of the first~$c-1$ edges of~$i_1\leftrightarrow j_1$. Let~$c'$ be the~$c$-th vertex of~$i_1\leftrightarrow j_1$.  By \Cref{cor:completionRCOP} there exists an automorphism~$\gamma\in \Gamma(\cG)$ that maps~$i_1\leftrightarrow c'$ to~$k_1\leftrightarrow c$. Substitute~$i_1\leftrightarrow j_1$ in~$L$ with~$k_1\leftrightarrow \gamma(j_1)$ by applying  the Markov move~$(i_1,j_1)-(\gamma(i_1),\gamma(j_1))$ 
from~$\mathcal{B}_{\overline{\cG}}$.  For simplicity, we update~$L$ so that the path~$i_1\lra j_1$ is replaced with~$k_1\lra \gamma(j_1)$.

If~$c$,~$l_1$ and the new~$j_1$ are all equal, the move has sent~$i_1\leftrightarrow j_1$ to~$k_1\leftrightarrow l_1$. 
Thus, the resulting binomial is
 $\sigma_{k_1 l_1} \big(\prod_{m=2}^d \sigma_{i_m j_m} - \prod_{m=2}^d \sigma_{k_m l_m} \big)$. As~$\ker(\varphi_{\cG})$ is toric, it is prime and contains no monomials. So,~$\prod_{m=2}^d \sigma_{i_m j_m} - \prod_{m=2}^d \sigma_{k_m l_m}$ is in~$ \ker(\varphi_{\cG})$, again. This binomial has degree~$d-1$. Thus,
by induction we know that it belongs to~$I_G + I_{\Gbar}$.

 Suppose now that at least one of~$j_1$ and~$l_1$ is strictly greater than~$c$. Choose that~$l_1> c$; this assumption is without loss of generality, since if we had $l_1 = c$, we could simply swap the roles of $L$ and $R$ and obtain the analogous results. If~$j_1$ is also greater than~$ c$,  then denote by~$a$ the vertex following~$c$ in~$c\lra j_1$. By our previous assumption {on~$c$}, we have~$\lambda(\{c,c+1\})\neq \lambda(\{c,a\})$ and that $c \lra j_1$ and $c \lra l_1$ intersect only at $c$. We are now precisely in the setting of \Cref{hyp1}.
By \Cref{lm:i2j2existance}, there cannot simultaneously be an edge in~$a\lra j_1$ of color~$\lambda(\{c,c+1\})$ and an edge in~$(c+1)\lra l_1$  of color~$\lambda(\{c,a\})$. 
Thus we may assume that the path~$a\leftrightarrow j_1$ does not contain an edge of color~$\lambda(\{c,c+1\})$; note that otherwise, exchanging the roles of~$i_1\leftrightarrow j_1$ and~$k_1\leftrightarrow l_1$ places us in an analogous situation.   
Since~$\lambda(\{c,c+1\})$ appears in~$\Lambda(c\lra l_1)$ but not in~$\Lambda(c\lra j_1)$, there must be another path in~$L$, say~$i_2\leftrightarrow j_2$, that uses the edge color~$\lambda(\{c,c+1\})$.  Since~$\cG$ is an RCOP graph, there is an automorphism~$\gamma$ in~$\Gamma(\cG)$ that maps this particular edge in~$i_2\leftrightarrow j_2$ to the edge~$\{c,c+1\}$. Substitute~$i_2\leftrightarrow j_2$ in~$L$ with~$\gamma(i_2)\leftrightarrow \gamma(j_2)$ using the linear Markov move~$(i_2,j_2)-(\gamma(i_2),\gamma(j_2))$ from~$\mathcal{B}_{\overline{\cG}}$. To simplify notation, we now refer to~$\gamma(i_2)$ and~$\gamma(j_2)$ as~$i_2$ and~$j_2$, respectively. If~$i_2 \neq c$, denote by~$b$ the vertex appearing before~$c$ in~$i_2\lra j_2$. We can now apply the previous lemmas of this section that rely on \Cref{hyp2}.

If $a$ or $b$ do not exist, or they exist but are non-adjacent, then \Cref{lm:nonadjacent} gives a Markov move that swaps $(i_1,j_1)$ with $(i_1,j_2)$.
If $a$ or $b$ are adjacent, then since they are also adjacent to $c$, $a,b$ and $c$ all belong to a maximal clique $\cC$. Either there is a another clique $\cC'$ that is isomorphic to $\cC$ and contains $c$, or there is not.
By \Cref{lem:severalcliques} and \Cref{lem:oneclique}, in each of these cases there exists a sequence of Markov moves in $\Bcal_G \cup \Bcal_{\overline{\Gcal}}$ 
that we may apply to $(i_1,j_1)\dots(i_d,j_d)$ in order to replace $(i_1,j_1)$ with $(i_1,j_1')$ such that $i_1 \lra j_1'$ passes through $c+1$.
Since~$i_1 = k_1$, we may repeat this procedure until~$k_1\lra l_1$ is in $L$.
Since~$k_1 \lra l_1$ belongs to both~$L$ and~$R$, and~$\ker(\varphi_{\cG})$ is prime and generated by binomials, we may now consider the binomial obtained by removing~$k_1 \lra l_1$ from both~$L$ and~$R$.
This binomial has degree~$d-1$, so we are done by induction.
\end{proof}

Now we are ready to finalize the proof of our main result. 

\begin{theorem}\label{thm:MainSec6}
Let~$\Gcal$ be an RCOP block graph with underlying graph~$G$ and completion $\hat{\cG}$. Then~$I_{\Gcal}=I_G+I_{\hat{\cG}}$.
In particular, $I_{\Gcal}$ is the toric ideal generated by the degree two binomials in~(\ref{eq:gensG}) and linear binomials in~(\ref{eq:completeRCOP}). 
\end{theorem}

\begin{proof}
\Cref{th:mainTheoremFormulation} shows that~$\Bcal_G \cup \Bcal_{\Gbar}$ forms a Markov basis for~$\ker(\varphi_{\cG})$.
Hence it suffices to show that~$I_{\Gcal} = \ker(\varphi_{\cG})$.

The linear space~$\Lcal_{\Gcal}$ can be obtained from either~$\Lcal_G$ or~$\Lcal_{\Gbar}$ by placing linear constraints on the parameters.
Hence,~$\Lcal_{\Gcal} \subset \Lcal_G, \Lcal_{\Gbar}$.
Inverting each linear space and taking the vanishing ideals yields that~$I_G, I_{\Gbar} \subset I_{\Gcal}$.
Thus the prime ideal~$\ker(\varphi_{\cG}) = I_G + I_{\Gbar}$ is contained in~$I_{\cG}$.
The ideals~$\ker(\varphi_{\cG})$ and~$I_{\Gcal}$ are prime as they are kernels of rational maps. So it suffices to show that they have the same dimension. 
In both ideals the number of colors used in~$\cG$ is an upper bound on the dimension as this is the number of parameters used in each parametrization. In the case of~$I_{\Gcal}$ this upper bound is tight as the spaces ~$\Lcal_\cG$ and~$\Lcal^{-1}_\cG$ have the same~dimension. Since~$\ker(\varphi_{\cG}) \subset I_{\Gcal}$ we have that $\dim(\ker(\varphi_{\cG})) \geq \dim(I_{\Gcal})$.  Hence they must have the same dimension. 
\end{proof}

\subsection{Discussion}\label{sec:discussion}
RCOP models belong to the larger class of \emph{colored Gaussian graphical models}. These are linear concentration models that arise from any colored graph $\Gcal$ wherein non-edge in $\Gcal$ correspond to zeros in $K$ and two entries of $K$ are set equal if and only their corresponding features have the same color in $\Gcal$.
The colored graph in   \Cref{ex:Frets_heads} suggests that disrupting block graph conditions in \Cref{thm:MainSec6} produces ideals~$I_{\Gcal}$ that are not toric. Moreover, in \Cref{ex:changeofcoords}, we will see that colored block graphs that are not RCOP can also have non-toric vanishing ideals. In fact, this is the case for every example that the authors computed. So we conjecture that toric vanishing ideals  arise only from RCOP block graphs.

\begin{conjecture}
\label{conj:toricCGM}
Let~$\Gcal$ be a colored graph. Then~$I_{\Gcal}$  is toric if and only if~$\Gcal$ is an RCOP block graph.
\end{conjecture}

The set~$\mathcal{L}^{-1}_{\cG}$, however, may still have a hidden toric structure even if its vanishing ideal~$I_{\cG}$ is not toric. To reveal this toric structure one needs to search for a different parametrization of~$\mathcal{L}^{-1}_{\cG}$. This can be done by applying a linear change of coordinates to~$\mathbb{R}[\Sigma]$ under which the ideal~$I_{\cG}$ becomes toric, as in the following example. 
\begin{example}\label{ex:changeofcoords}
Let~$\cG$ be the colored path on three vertices where each vertex has a distinct color and where the two edges~$[1,2]$ and~$[2,3]$ share the same color. 
Its  vanishing ideal is
\(
I_{\cG}= \langle \sigma_{13}\sigma_{22}-\sigma_{12}\sigma_{23},\sigma_{12}\sigma_{13}-\sigma_{11}\sigma_{23}-\sigma_{13}\sigma_{23}+\sigma_{12}\sigma_{33} \rangle
\)
which is clearly not toric  However, its image under the  linear change of coordinates~$p_{11}= \sigma_{11}+\sigma_{13}, \ p_{33}= \sigma_{33}+\sigma_{13}$, and~$p_{ij}= -\sigma_{ij}$ for all other $i,j\leq 3$, 
is the toric ideal~$J = \langle p_{13}p_{22}-p_{12}p_{23}, p_{11}p_{23}-p_{12}p_{33} \rangle$.
Thus the variety $V(I_{\Gcal})$ is equal to the image of $V(J)$ under an invertible linear transformation.
\end{example}

This ends the discussion with the hopeful note that  the class of colored Gaussian graphical models with underlying toric structure extends outside the RCOP block graphs and is worth exploring.

\bibliographystyle{plain}
\bibliography{paper}

\hfill

\noindent
\footnotesize {\bf Authors' addresses:}

\hfill

\noindent Jane Ivy Coons\\
St John's College, University of Oxford, United Kingdom\\
\hfill {\tt jane.coons@maths.ox.ac.uk}

\hfill

\noindent Aida Maraj\\
Department of Mathematics, \\ Universiry of Michigan, Ann Arbor (MI), United States of America\\
\hfill {\tt maraja@umich.edu}

\hfill

\noindent Pratik Misra\\
Department of Mathematics, \\ KTH Royal Institute of Technology, SE-100 44 Stockholm, Sweden \\
\hfill {\tt pratikm@kth.se}

\hfill

\noindent Miruna-\c Stefana Sorea\\ SISSA - Scuola Internazionale Superiore di Studi Avanzati, Trieste, Italy and RCMA Lucian Blaga University, Sibiu, Romania\\
\hfill {\tt mirunastefana.sorea@sissa.it}
\end{document}